\documentclass[12pt]{article}
\usepackage{amsmath}
\usepackage{amsfonts}
\usepackage{amssymb}
\usepackage{polski}

\input{tcilatex}
\begin{document}

\begin{center}
\textbf{The Goldie Equation: III. Homomorphisms from Functional Equations.}

\textbf{by}

\textbf{N. H. Bingham and A. J. Ostaszewski.}

\bigskip
\end{center}

\textbf{Abstract}. This is the second of three sequels to [Ost 3] -- the
third of the resulting quartet -- concerning the real-valued continuous
solutions of the multivariate Goldie functional equation $(GFE)$ below of
Levi-Civita type. Following on from the preceding paper [BinO8], in which
these solutions are described explicitly, here we characterize $(GFE)$ as
expressing homomorphy (in all but some exceptional \textquotedblleft
improper\textquotedblright\ cases) between multivariate Popa groups, defined
and characterized earlier in the sequence. The group operation involves a
form of affine addition (with local scalar acceleration) similar to the
circle operation of ring theory. We show the affine action in $(GFE)\ $may
be replaced by general continuous acceleration yielding a functional
equation $(GGE)$ which it emerges has the same solution structure as $(GFE)$%
. The final member of the sequence [BinO9] considers the richer framework of
a Banach algebra which allows vectorial acceleration, giving the closest
possible similarity to the circle operation.

\bigskip

\textbf{1. Introduction}. We begin by recalling the classic Go\l \k{a}%
b-Schinzel equation 
\begin{equation}
\eta (x+\eta (x)y)=\eta (x)\eta (y)\qquad (x,y\in X)  \tag{$GS$}
\end{equation}%
for $X=\mathbb{R}.$ Interpreted in the context of a real topological vector
space $X,$ the continuous solutions of $(GS)$ are in one of the following
two forms, both generalizing those for $X=\mathbb{R}$. One is%
\begin{equation*}
\eta (x)=1+\rho (x)\qquad (x\in X),
\end{equation*}%
for some continuous linear functional $\rho :X\rightarrow \mathbb{R}$. The
other is%
\begin{equation*}
\eta (x)=\max \{1+\rho (x),0\}\qquad (x\in X).
\end{equation*}%
(This is the Brillou\"{e}t-Dhombres-Brzd\k{e}k theorem [BriD, Prop. 3],
[Brz1, Th. 4]; cf. [Brz2]) Significantly for us, the latter form is
non-negative. We assume henceforth that $\rho \neq 0.$

Next, recall from algebra the \textit{circle operation} of ring theory (see
e.g. Jacobson [Jac1951], II.3), which in a ring identifies the elements
which have inverses, and so form a \textit{group}:%
\begin{equation*}
a\circ b:=a+b-ab.
\end{equation*}

Analogously, a \textit{Popa group} on a real topological vector space $X$ is
obtained from a given $\rho \in X^{\ast }$, a continuous linear map $\rho
:X\rightarrow \mathbb{R}$, by defining%
\begin{equation*}
u\circ _{\rho }v:=u+v+\rho (u)v=u+\eta _{\rho }(u)v,
\end{equation*}%
where%
\begin{equation}
\eta _{\rho }(u):=1+\rho (u).  \tag{$\eta _{\rho }$}
\end{equation}%
This equips%
\begin{equation*}
\mathbb{G}^{+}(X)=\mathbb{G}_{\rho }(X):=\{x\in X:\eta _{\rho }(x)>0\}
\end{equation*}%
with a group structure. (Below the notation $\mathbb{G}^{+}(X)$ is useful
whenever context dictates that $\eta _{\rho }$ or some similar function,
such as $h$ below, is to be positive.) Its significance for $(GS)$ was first
recognized by Popa [Pop]. Its explanatory power as a fundamental tool in the
study of regular variation is witnessed in several papers (both joint
[BinO2-7] and separate [Ost1-4]) and especially in combination with
shift-compactness -- see [BinO10] and the survey [BinO11]; cf. [BinO1].

We recall the associated \textit{Goldie functional equation}, originally of
(univariate) regular variation [BinO3] for the pair $(K,g)$:%
\begin{equation}
K(u\circ _{\rho }v)=K(u)+g(u)K(v)\qquad (u,v\in \mathbb{G}_{\rho }(X)), 
\tag{$GFE$-$\rho $}
\end{equation}%
where now with $g:X\rightarrow \mathbb{R}_{+}:=[0,\infty ),$ $K:\mathbb{G}%
_{\rho }(X)\rightarrow Y,$ for $Y$ again a real topological vector space,
and with both functions \textit{continuous}.

Given the pair, $K$ is called the \textit{kernel, }$g$ the (outer) \textit{%
auxiliary }function of the kernel, following usage in regular variation
theory, where this equation plays a key role, see e.g. [BinO8]. It is
helpful to view $g$ as an \textit{action} imparting local \textit{%
acceleration} to the action of addition.\textit{\ }The equation is a special
case of the \textit{Levi-Civita functional equation} [Lev], [Stet, Ch. 5]
(for the wider literature, see [AczD, Ch. 14, 15]), usually studied on
semi-groups. We shall also study the \textit{generalized Goldie equation}
for the triple $(K,h,g)$, now with inner and outer auxiliaries:%
\begin{equation}
K(u+h(u)v)=K(u)+g(u)K(v)\qquad (u,v\in \mathbb{G}^{+}(X)),  \tag{$GGE$}
\end{equation}%
with $g,h:X\rightarrow \lbrack 0,\infty ),$ $K:X\rightarrow Y$ and all three 
\textit{continuous }with%
\begin{equation*}
\mathbb{G}^{+}(X)=\mathbb{G}_{h}^{+}(X):=\{x\in X:h(x)>0\},
\end{equation*}%
using analogous notation. It emerges in Lemma 8.1 that $h$ here may be
replaced by $g$ (i.e. briefly, that $\mathbb{G}_{h}^{+}(X)=\mathbb{G}%
_{g}^{+}(X)).$

When $\rho \in X^{\ast }$ is fixed, $(GGE)$ reduces to $(GFE)$ for $u,v\in 
\mathbb{G}_{\rho }(X)$ when%
\begin{equation*}
h(u)=\eta _{\rho }(u)=1+\rho (u).
\end{equation*}

For fixed $\rho \in X^{\ast },K:\ \mathbb{G}_{\rho }(X)\rightarrow Y$ and a
given continuous linear map $\sigma :Y\rightarrow \mathbb{R}$, take $%
g(u)=g^{\sigma }(u),$ for $u\in \mathbb{G}_{\rho }(X),$ where 
\begin{equation}
g^{\sigma }(u):=\eta _{\sigma }(K(u))=1+\sigma (K(u)).  \tag{$g^{\sigma }$}
\end{equation}%
Then $(GFE)$ for $(K,g)$ reduces to the Popa \textit{homomorphism} [BinO5]
from $\mathbb{G}_{\rho }(X)$ to $\mathbb{G}_{\sigma }(Y):$%
\begin{equation*}
K(u\circ _{\rho }v)=K(u)\circ _{\sigma }K(v)\quad (=K(u)+K(v)+\sigma
(K(u))K(v)).
\end{equation*}%
Theorem 7.2 asserts that if $(K,g)$ satisfies $(GFE),$ then $g=g^{\sigma }$
for some $\sigma ,$ unless the range $\mathcal{R}(K)$ collapses to $K(%
\mathcal{N}(\rho )).$

\bigskip

\textbf{Notational convention:} Throughout below, $X,Y$ are fixed real
topological vector spaces; $\rho \in X^{\ast }$ may either be given in
advance, or specifically constructed. For given $\rho \in X^{\ast },$ we say
that the pair $(K,g)$ \textit{satisfies} $(GFE)$ to mean that $%
g:X\rightarrow \mathbb{R}_{+}:=(0,\infty ),$ $K:\mathbb{G}_{\rho
}(X)\rightarrow Y,$ with both functions \textit{continuous,} and that $(GFE)$
is satisfied by the pair $(K,g)$. Always, `$(GFE)$' is to be read as short
for `$(GFE$-$\rho )$'; both imply the presence of a given $\rho \in X^{\ast
} $ and the corresponding domain-restriction to $\mathbb{G}_{\rho }(X).$

Likewise, we say that the triple $(K,h,g)$ \textit{satisfies} $(GGE),$ to
mean that $g,h:X\rightarrow \mathbb{R}_{+},$ $K:X\rightarrow Y,$ with all
three functions \textit{continuous} and that the relation $(GGE)$ is
satisfied by the triple $(K,h,g)$. In Section 8 it is shown that this
implies that the inner auxiliary $h$ is given by $(\eta _{\rho })$ above for
some $\rho $.

\bigskip

For (fixed $\rho \in X^{\ast }$ and) a given pair $(K,g)$ satisfying $(GFE)$
the existence theorem, Theorem 7.2 below, asserts that, unless the range $%
\mathcal{R}(K)$ collapses to $K(\mathcal{N}(\rho )),$ there \textit{exists}
a \textit{unique} linear $\sigma =\sigma _{g}:Y\rightarrow \mathbb{R}$ with $%
g=g^{\sigma }$ as in ($g^{\sigma }$) above, which is \textit{continuous}
provided the range $K(\mathbb{G}_{\rho }(X))$ is a \textit{closed
complemented subspace} of $Y$ (see \S 7). The argument is involved and
begins by establishing necessary and sufficient conditions on any such $%
g:X\rightarrow \mathbb{R}_{+}$ that there exists a linear $\sigma
:Y\rightarrow \mathbb{R}$ with $g=g^{\sigma }$ as above. In an intermediate
step (see Prop. 7.1A and 7.1B) we deduce that, unless the range $\mathcal{R}%
(K)\ $\textit{collapses} to $K(\mathcal{N}(\rho )),$ such a $\sigma $ always
exists and is unique, so that we may refer to it as $\sigma _{g},$ and that
this $\sigma _{g}$ is continuous provided $K(\mathbb{G}_{\rho }(X))$ is a
closed complemented subspace of $Y.$ The result thus extends to continuous
functions on arbitrary linear domains the one-dimensional homomorphy first
recognized in [Ost3, Th. 1].

It emerges that intersections of the null spaces of various additive maps on 
$X$ are of central significance. So, given $g:X\rightarrow \mathbb{R}_{+},$
we put for $x\in \mathbb{G}_{g}^{+}(X)$%
\begin{equation*}
\gamma (x):=\log g(x),
\end{equation*}%
and for additive $\alpha :\mathbb{G}_{\rho }(X)\rightarrow \mathbb{R}$ we
write 
\begin{equation*}
\mathcal{N}(\alpha ):=\{x\in X:\alpha (x)=0\},\qquad \mathcal{N}^{\ast
}(\alpha ):=\mathcal{N}(\alpha )\cap \mathcal{N}(\rho );
\end{equation*}%
of particular interest here is%
\begin{equation*}
\mathcal{N}^{\ast }(\gamma ):=\mathcal{N}(\gamma )\cap \mathcal{N}(\rho
)=\{x:g(x)=1\}\cap \mathcal{N}(\rho ).
\end{equation*}%
So we begin in Section 2 with a study of the\textit{\ }auxiliary function $g$
corresponding to a given pair $(K,g)$ satisfying $(GFE),$ leading us to
Theorem 2.1, which we view as establishing an \textit{index} (cf. [BinO7]).
Section 3, prompted by the radial properties of Popa homomorphisms
established in [BinO8], asserts in Theorem 3.1 analogous radiality
properties of $(GFE)$ kernels. The proof is delayed to Section 5, after
establishing in Section 4 the density of sets canonically modelling the
rationals $m/n$ by appropriate `steering' of $m$-fold $g$- or $h$-actions
applied to $u/n$ and computing related limits. This radiality implies
(Corollary 6.2) that any kernel $K:\mathbb{G}_{\rho }(X)\rightarrow Y$
satisfies a dichotomy involving the null spaces $\mathcal{N}(\rho )$ and $%
\mathcal{N}(\alpha )$ for $\alpha :\mathbb{G}_{\rho }(X)\rightarrow \mathbb{R%
}$ an arbitrary additive map. The dichotomy concerning the two null spaces
arises because two hyperplanes passing through the origin (representing the
two null spaces) have intersection either with co-dimension 1, when they
coincide, or 2 otherwise.

\bigskip

In Section 7 we induce a \textit{Popa group structure} on the range space $Y$
and prove our main result, Theorem 7.2, on the existence of an appropriate
functional $\sigma \in Y^{\ast }$ for transforming the right-hand side of $%
(GFE)$ into the circle operation associated with $\mathbb{G}_{\sigma }(Y)$.
Section 8 is devoted to establishing a reduction to $(GFE)$ of the more
general Goldie functional equation $(GGE)$ above, wherein the accelerated
summation $u+h(u)v$ on the left-hand side of $(GGE)$ replaces the Popa
operation $u\circ _{\rho }v$ . The argument of Section 8 relies on radial
behaviours and on reducing a known `pexiderised variant' \footnote{%
That is, additional function symbols replace instances in the functional
equation $(GS)$ of the function symbol of primary interest.} of $(GS)$ to $%
(GS)$ itself, with a resulting \textit{tetrachotomy} of possible Popa
homomorphisms in any direction $u\in X$; the latter foursome can be
interpreted as arising from the binary split into the vanishing or
non-vanishing of Gateaux derivatives (along the radial direction $u$) of $h$
and $g,$ whose existence follows from $h$ satisfying $(GS)$ and $g$
satisfying the pexiderized equation $(PGS)$ (see \S 8).

\bigskip\ 

\noindent \textbf{2. Auxiliary functions: multiplicative property. }We learn
from Lemma 2.1 that the auxiliary function $g$ of a kernel is $\rho $%
-multiplicative in the sense of $(M)$ below, and so the corresponding $%
\gamma =\log g$ (see above) is $\rho $-additive. When context permits we
omit the prefix.

\bigskip

\noindent \textbf{Lemma 2.1 }(cf. [BinO5], [Ste, Prop. 5.8])\textbf{. }%
\textit{If }$(K,g)$\textit{\ satisfies }$(GFE)$\textit{\ and }$K$\textit{\
is non-zero, then }$g$ \textit{is }$\rho $\textit{-multiplicative:}%
\begin{equation}
g(u\circ _{\rho }v)=g(u)g(v)\qquad (u,v\in \mathbb{G}_{\rho }(X)),  \tag{$M$}
\end{equation}%
\textit{and so }$\gamma =\log g$\textit{\ is }$\rho $\textit{-additive:}%
\begin{equation}
\gamma (u\circ _{\rho }v)=\gamma (u)+\gamma (v)\qquad (u,v\in \mathbb{G}%
_{\rho }(X)).  \tag{$A$}
\end{equation}%
\textit{Fix }$u\neq 0.$ \textit{For }$t\in \mathbb{R},$\textit{\ if }$%
g(tu)\neq 1$\textit{\ except at }$t=0,$\textit{\ then }$g(tu)$\textit{\
takes one of two forms:}%
\begin{equation*}
g(tu)=\left\{ 
\begin{array}{cc}
(1+t\rho (u))^{\gamma (u)/\log (1+\rho (u))}, & \rho (u)>0, \\ 
e^{\gamma (u)t}, & \rho (u)=0.%
\end{array}%
\right.
\end{equation*}

\bigskip

\noindent \textbf{Proof. }To deduce $(M)$ use the two ways to associate the
three variables in $K(u\circ _{\rho }v\circ _{\rho }w)$. For the routine
details see \S 9.1 (Appendix). Put $g_{u}(t):=g(tu);$ then $g_{u}:\mathbb{G}%
_{\rho (u)}\mathbb{(R)\rightarrow R}_{+}$ and satisfies $(M)$ on the real
line. This case is covered by established results, e.g. [BinO5] or Th. BO
below, yielding for some $\kappa (u)$%
\begin{equation*}
g_{u}(t)=(1+\rho (u)t)^{\kappa (u)/\rho (u)},
\end{equation*}%
whence the cited formulas. \hfill $\square $

\bigskip

\noindent \textbf{Remark. }In \S 3 below we encounter the \textit{link
function} $\lambda $, and with it ${\lambda }_{u}$, in terms of which we
will be able to write, in view of Theorem 3.1., 
\begin{equation}
K_{u}(tu):=\lambda _{u}(t)K(u).  \tag{$K$}
\end{equation}%
Then, for $K(u)\neq 0,$%
\begin{equation*}
\lambda _{u}(s\circ _{\rho }t)=\lambda _{u}(t)+g(tu)\lambda _{u}(t),
\end{equation*}%
whence, from the real context of [BinO5], $K_{u}(tu)$ is proportional to one
of the two forms%
\begin{equation*}
\left\{ 
\begin{array}{cc}
\lbrack (1+t\rho (u))^{\gamma (u)}-1], & \text{ if }\rho (u)>0, \\ 
\lbrack e^{\gamma (u)t}-1], & \text{ if }\rho (u)=0.%
\end{array}%
\right.
\end{equation*}%
This is implied by Lemma 6.1, itself a corollary of Theorem 2.1 below.

The radial case $g_{u}$ of Lemma 2.1 shows that we need not restrict to
scalars. We can further describe $g(x)$ explicitly by studying its
associated \textit{index} $\gamma (x)$, to borrow a term from extreme-value
theory (EVT), for which see e.g. [BinO7], [HaaF, p. 295]. Below we
distinguish between linearity (in the sense of vectors and scalars), which $%
\gamma $ exhibits \textit{only} on $\mathcal{N}(\rho ),$ and its more
general property of $\rho $-additivity. To see the difference note that, for
distinct $u$ and $w$ with $\rho (u)=1$ and $\rho (w)=1$ (so that $w-u\in 
\mathcal{N}(\rho )$), taking $x=tw$ in the first display in Theorem 2.1
below gives%
\begin{eqnarray*}
g(tw) &=&g(x)=e^{t\gamma (w-u)}(1+t)^{\gamma (u)/\log 2}=(1+t)^{\gamma
(w)/\log 2}: \\
\qquad e^{t\gamma (w-u)} &=&(1+t)^{[\gamma (w)-\gamma (u)]/\log 2}.
\end{eqnarray*}%
We thus think of the following result as an Index theorem.

\bigskip

\noindent \textbf{Theorem 2.1.} \textit{For }$(K,g)$ \textit{satisfying }$%
(GFE),$\textit{\ the auxiliary function }$g$\textit{\ is }$\rho $\textit{%
-multiplicative and its index }$\rho $\textit{-additive.}

\textit{So, as in Lemma 2.1, for any }$u$ \textit{with }$\rho (u)=1,$%
\begin{equation*}
g(x)=e^{\gamma (x-\rho (x)u)}(1+\rho (x))^{\gamma (u)/\log 2}\qquad (x\in 
\mathbb{G}_{\rho }(X)),
\end{equation*}%
\textit{where, for }$\gamma =\log g,$ 
\begin{equation*}
\alpha (x):=\gamma (x-\rho (x)u)
\end{equation*}%
\textit{\ \noindent is linear and }$\alpha (u)=0.$

\textit{Conversely, for any }$\alpha :\mathbb{G}_{\rho }(X)\rightarrow 
\mathbb{R}$\textit{\ additive and }$\beta $ \textit{a real parameter, the
following function is multiplicative (satisfies }$(M)$\textit{):}%
\begin{equation*}
\bar{g}(x)=\bar{g}_{\alpha ,\beta }(x):=e^{\alpha (x)}(1+\rho (x))^{\beta }.
\end{equation*}%
\noindent \textbf{Proof.} By Lemma 2.1, $\gamma $ satisfies ($A$). So $%
\gamma :\mathbb{G}_{\rho }(X)\rightarrow \mathbb{G}_{0}(\mathbb{R})=(\mathbb{%
R},+).$ Here $\gamma (\mathcal{N}(\rho ))\subseteq \mathbb{R}=\mathcal{N}%
(0). $ By a theorem of Chudziak [Chu3, Th. 1] as amended in [BinO8, Th. Ch.,
Th. 4A], for any $u$ with $\rho (u)=1$,%
\begin{eqnarray*}
\gamma (x) &=&\gamma (x-\rho (x)u)+[\gamma (u)/\log 2]\log [(1+\rho (x))]: \\
g(x) &=&e^{\gamma (x-\rho (x)u)}(1+\rho (x))^{\gamma (u)/\log 2}.
\end{eqnarray*}%
Then, taking $x=tu$ $,$%
\begin{equation*}
g(tu)=(1+t)^{\gamma (u)/\log 2},
\end{equation*}%
as, by linearity, $\rho (tu)=t.$ For $w$ with $\rho (w)>0,$ take $u=w/\rho
(w).$ Then 
\begin{equation*}
g(tw)=g(t\rho (w)u)=(1+t\rho (w))^{\gamma (u)/\log 2}.
\end{equation*}%
On the other hand, for $w$ with $\rho (w)=0,$%
\begin{equation*}
g(tw)=e^{\gamma (w)t}.
\end{equation*}%
Given the form of $\bar{g},$ it is routine to check that $(M)$ holds. See \S %
9.2 (Appendix) for details.\hfill $\square $

\bigskip

\noindent \textbf{Lemma 2.2.} \textit{For continuous }$(K,g)$\textit{\
satisfying }$(GFE)$\textit{\ and }$\gamma =\log g,$\textit{\ }$\mathcal{N}%
(\gamma )$\textit{\ is a subgroup of }$\mathbb{G}_{\rho }(X),$\textit{\ and }%
$\mathcal{N}^{\ast }(\gamma )$\textit{\ is both a vector subspace and a }$%
\mathbb{G}_{\rho }(X)$\textit{-subgroup of }$\mathcal{N}(\rho )$.

\bigskip

\noindent \textbf{Proof. }By Lemma 2.1, $\gamma $ is additive on $\mathbb{G}%
_{\rho }(X)$ and so linear on $\mathcal{N}(\gamma ),$ by continuity of $%
\gamma .$ The remaining assertions are clear. \hfill $\square $

\bigskip

\noindent \textbf{3. Radiality. }We write $\langle v\rangle $ for the linear
span of a vector $v$, in $X$ or $Y$. For $u\in \mathbb{G}_{\rho }(X),$ we
write%
\begin{equation*}
\langle u\rangle _{\rho }:=\langle u\rangle \cap \mathbb{G}_{\rho
}(X)=\{tu:t\in \mathbb{R},1+\rho (tu)>0\}.
\end{equation*}%
Our main result here and our later tool is Theorem 3.1 below. This asserts 
\textit{radiality}, the property that the kernel function maps the points
along $\langle u\rangle $ to points along $\langle K(u)\rangle ,$ and,
furthermore, specifies precisely the \textit{linkage} between the
originating vector $tu$ and its image $K(tu)=\lambda _{u}(t)K(u)$, as in $%
(K).$ The dependence is uniform, through one and the same \textit{link
function} $\lambda $ (below) but with its continuously varying parameters
referring to what we term informally the \textit{growth rates} (below) of
the two auxiliaries along $u$ at the origin. We use either the notation $%
\lambda (t;r,\theta )$ which identifies the parameters explicitly, or $%
\lambda _{u}$ when the parameters are implied by the direction $u.$

To state it in a form adequate to cover both $(GFE)$ and $(GGE),$ we need
several definitions and a lemma. Notice that $(GGE)$ implies that 
\begin{equation*}
K(h(0)v)=K(0)+g(0)K(v),
\end{equation*}%
so if $h(0)=0$ and $g(0)\neq 0,$ the map $K\ $is trivial; similarly, if $%
g(0)=0$ and $h(0)\neq 0,$ since Theorem 3.1 below asserts that $%
K(tv)=\lambda (t)K(v)$ for some monotone function $\lambda .$ (strictly
monotone if $K(v)\neq 0).$ Thus without loss of generality (briefly,
w.l.o.g.) we standardize the auxiliary functions in $(GGE),\ $by taking $%
h(0)=g(0)=1.$ This then coincides with the corresponding conditions for $%
(GFE)$ in Theorem 2.1. It now easily follows that $K(0)=0,$ since%
\begin{equation*}
K(0)=K(0+h(0)0)=K(0)+g(0)K(0).
\end{equation*}

\noindent \textbf{Definition.} 1. Following Th. 2.1 above, taking as
parameters $r\geq 0,\theta \in \mathbb{R}$, the standard \textit{%
multiplicative radial function} $g_{r,\theta }(t)$ is defined for $t>-1/r$
(the latter interpreted for $r=0$ as $-\infty )$ by%
\begin{equation*}
g_{r,\theta }(t):=\left\{ 
\begin{array}{cc}
(1+tr)^{\theta /r}, & \text{if }r>0,\text{ } \\ 
e^{\theta t}, & \text{if }r=0.%
\end{array}%
\right.
\end{equation*}%
Thus $g_{r,\theta }$ is Gateaux differentiable at $t=0$. This definition
blends two possible instances of the function $\bar{g}$ of Theorem 2.1
consistently with the L'Hospital convention (as $g_{0,\theta
}=\lim_{r\searrow 0}g_{r,\theta }),$ constantly applied below implicitly.

\noindent 2. Taking again as parameters $r\geq 0$ and $\theta \in \mathbb{R}$%
, we define below the function $\lambda (t;r,\theta )$ for $t>-1/r$ (where
for $r=0,$ we interpret $-1/r$ as $-\infty $), which we call the \textit{%
Popa link }function on account of its role in Theorem 3.1 below:%
\begin{equation}
\lambda (t)=\lambda (t;r,\theta ):=\left\{ 
\begin{array}{cc}
\lbrack (1+rt)^{\theta /r}-1]/[(1+r)^{\theta /r}-1], & \text{if }r>0,\theta
\neq 0, \\ 
\ln (1+rt)/\ln (1+r), & \text{if }r>0,\theta =0, \\ 
(e^{t\theta }-1)/(e^{\theta }-1), & \text{if }r=0,\theta \neq 0, \\ 
t, & \text{if }r=\theta \in \{0,\pm \infty \}.%
\end{array}%
\right.  \tag{$\dag $}
\end{equation}

This function first arises in the context of $GFE$ as a map $\mathbb{G}_{r}(%
\mathbb{R})\rightarrow \mathbb{R}$ with parameters $r,$ $\theta $ and $%
\lambda $ satisfies an analogous equation $(GFE_{\lambda })$ below. But as
it plays an equivalent role in $(GGE),$ it is convenient to derive its
properties in the more general setting. Its key elementary properties are
summarized in the following lemma; see also Lemma 8.4.

\bigskip

\noindent \textbf{Lemma 3.1.} \textit{For} $r\in \lbrack 0,\infty ),\theta
\in \mathbb{R}$, \textit{the Popa link function }$\lambda $\textit{\ is
separately continuous in }$t$\textit{\ and in its parameters, with }$\lambda
(0)=0$\textit{\ and }$\lambda (1)=1,$\textit{\ and satisfies the equation}%
\begin{equation}
\lambda (s\circ _{r}t)=\lambda (s)+g_{r,\theta }(s)\lambda (t)\qquad (s,t\in 
\mathbb{G}_{r}(\mathbb{R})),  \tag{$GFE_{\lambda }$}
\end{equation}%
\textit{\ equivalently}%
\begin{equation*}
\lambda (s\circ _{r}t)=\lambda (s)\circ _{\sigma }\lambda (t),\text{ for }%
\sigma =g_{r,\theta }(1)-1.
\end{equation*}%
\textit{Except for }$r=\theta =0,$\textit{\ the equation }$\lambda (t)=t$%
\textit{\ has a unique solution }$t=1.$

\bigskip

\noindent \textbf{Proof.} Routine: see \S 9.3 (Appendix). \hfill $\square $

\bigskip

The main result of this section is Theorem 3.1 below, asserting the
radiality property of the kernel function in $(GFE)$ and $(GGE)\ $that, for
some scalar function $\lambda =\lambda _{u}$ which we determine,%
\begin{equation*}
K(su)=\lambda (s)K(u)\text{ for }s\geq 0.
\end{equation*}%
It is convenient to use notation bringing together some arguments common to $%
g$ and $h.$ Put%
\begin{eqnarray*}
\delta _{n}(u) &:&=\delta _{n}^{g}=g(u/n)-1,\text{ or }\delta
_{n}^{h}=h(u/n)-1, \\
\gamma _{g}(u) &:&=\lim\nolimits_{n}n\delta _{n}^{g}(u),\qquad \text{ and }%
\gamma _{h}(u)=\lim\nolimits_{n}n\delta _{n}^{g}(u),
\end{eqnarray*}%
whenever these limits exists, possibly $\pm \infty .$

In the context $h(u)=1+\rho (u)$ of $(GFE),$ we obtain $n\delta
_{n}^{h}=\rho (u),$ so that $\gamma _{h}(u)=\rho (u).$ Recall from \S 2 that 
$\gamma =\log g;$ here $g(su)=e^{\gamma (su)}$ gives $\gamma _{g}(u)=$ $%
\gamma (u),$ motivating the common notation.

\bigskip

\noindent \textbf{Theorem 3.1. }\textit{If }$(K,h,g)$\textit{\ satisfies }$%
(GGE)$\textit{\ with }$K,g,h$\textit{\ continuous and}%
\begin{equation*}
n\delta _{n}^{g}(u)\rightarrow \gamma _{g}(u)\in \mathbb{R}\text{ and }%
n\delta _{n}^{h}\rightarrow \gamma _{h}(u)\in \mathbb{R},
\end{equation*}%
\textit{then for }$u$ \textit{with} $K(u)\neq 0$%
\begin{equation*}
K\left( \frac{e^{\gamma _{h}(u)t}-1}{\gamma _{h}(u)}u\right) =\frac{%
e^{\gamma _{g}(u)t}-1}{e^{\gamma _{g}(u)}-1}K\left( \frac{e^{\gamma
_{g}(u)}-1}{e^{\gamma _{g}(u)}}u\right) .
\end{equation*}
\textit{In particular, if }$\gamma _{g}(u)=\gamma _{h}(u),$\textit{\ then }%
\begin{equation*}
K(tu)=tK(u)\text{ }\qquad (t\in \mathbb{R}).
\end{equation*}%
\textit{This applies also if one or other limit is infinity, both then being
equal. Furthermore, the usual L'Hospital convention applies when either of
the limits }$\gamma _{g}(u),\gamma _{h}(u)$ \textit{is zero. Thus}%
\begin{equation}
K(tu)=\lambda _{u}(t)K(u)\qquad \text{ for }\lambda _{u}(t):=\lambda
(t;\gamma _{h}(u),\gamma _{g}(u)).  \tag{$Rad$}
\end{equation}%
\textit{In particular, }$K$\textit{\ is Gateaux differentiable in direction }%
$u$\textit{: }%
\begin{equation*}
K^{\prime }(tu)=\lambda _{u}^{\prime }(t)K(u),
\end{equation*}%
\textit{and the function }$\lambda _{u}$ \textit{satisfies a corresponding
scalar }$(GGE)$:%
\begin{equation*}
\lambda _{u}(s+h(su)t)=\lambda _{u}(s)+g_{\gamma _{g}(u),\gamma
_{h}(u)}(su)\lambda _{u}(t).
\end{equation*}

The proof, based on Lemma 4.1 below, must be delayed until after a series of
limit calculations ending in Prop. 4.3. These rest only on the assumption
that $g$ and $h$ are \textit{continuous}. A stronger assumption (invoking 
\textit{Gateaux differentiability} of $g$ and $h$ at $0$) allows a much
shorter proof, directly from Lemma 4.1; see \S 9.4 (Appendix).

We will refer informally to the limits $\gamma _{g}(u)$ and $\gamma _{h}(u)$
above as the \textit{growth rates} of $g$ and $h.$

\bigskip

\noindent \textbf{4. Density preliminaries for Theorem 3.1}

The general approach of using a dense set to identify an unknown function is
familiar (see [Acz, Ch. 2, Ch. 6]), although here it is more involved, in
view of only a latent group structure (cf. [Ste]). Our first step is Lemma
4.1 for which we need a definition. This is phrased with a view to
generalizations beyond real-valued functions, allowing for $(GFE)$ to be
interpreted also in a Banach algebra setting.

\bigskip

\noindent \textbf{Definition. }The polynomials and rational polynomials (in
the indeterminate $x)$ $\wp _{n}$ and $[\wp _{m}/\wp _{n}]$ for $m,n\in 
\mathbb{N}$ are defined by:%
\begin{equation*}
\wp _{m}(x):=1+x+...+x^{m-1},\qquad \lbrack \wp _{m}/\wp _{n}](x):=\wp
_{m}(x)/\wp _{n}(x).
\end{equation*}%
So $\wp _{m}(1)=m,$ and so $[\wp _{m}/\wp _{n}](1)=m/n,$ and also $\wp
_{m}(t)=(t^{m}-1)(t-1)^{-1}$ when $t-1$ is invertible$.$

\bigskip

\noindent \textbf{Lemma 4.1}. \textit{If }$(K,h,g)$\textit{\ satisfies }$%
(GGE)$\textit{, then for }$u$ \textit{with} $K(u)\neq 0:$%
\begin{eqnarray}
K(\wp _{m}(h(u/n))u/n) &=&\wp _{m}(g(u/n))K(u/n),  \TCItag{$\ast $} \\
K(\wp _{m}(h(u/n))u/n) &=&[\wp _{m}/\wp _{n}](g(u/n))K(\wp _{n}(h(u/n))u/n).
\TCItag{$\ast \ast $}
\end{eqnarray}

\noindent \textbf{Proof. }For $x,y\in X,$ and $z$ an indeterminate ranging
over either $\mathbb{R}$ (as here) or a unital commutative Banach algebra $%
\mathbb{A}$ (as in [BOst9]), put%
\begin{equation*}
x\circ _{z}y:=x+yz.
\end{equation*}%
Starting from $u$ and $v:=K(u),$ we define a pair of sequences of `powers',
by iterating the operation $\circ _{z}$ for respectively $z=h(u)$ and $%
z=g(u).$ These iterates are defined inductively: 
\begin{equation*}
u_{h}^{n+1}=u\circ _{h(u)}u_{h}^{n},\qquad v_{g}^{n+1}=v\circ
_{g(u)}v_{g}^{n},\qquad \text{with }u_{h}^{1}=u,\qquad v_{g}^{1}=v.
\end{equation*}%
Then, for $n\geq 1,$%
\begin{equation}
K(u_{h}^{n+1})=K(u)+g(u)K(u_{h}^{n})=K(u)_{g}^{n+1}=v_{g}^{n+1}. 
\tag{$\ast
\ast \ast $}
\end{equation}%
Motivated by the case%
\begin{equation*}
K(u_{h}^{2})=K([1+h(u)]u)=K(u)+g(u)K(u)=[1+g(u)]K(u),
\end{equation*}%
the recurrence $(\ast \ast \ast )$ justifies associating with the iterates
above sequences of `coefficients' $(g_{n}(.)),$ $(h_{n}(.))$, by writing 
\begin{equation*}
v_{g}^{n}=g_{n}(u)K(u),\qquad u_{h}^{n}=h_{n}(u)u:\qquad
K(h_{n}(u)u)=g_{n}(u)K(u).
\end{equation*}

Solving appropriate recurrences arising from $(\ast \ast \ast )$ for the
iterations $u_{h}^{n+1}=u\circ _{h}u_{h}^{n}$ and $v_{g}^{n+1}=v\circ
_{g}v_{g}^{n}$ gives%
\begin{equation*}
u_{h}^{n}=h_{n}(u)u,\text{ and }v_{g}^{n}=g_{n}(u)K(u),
\end{equation*}%
where%
\begin{eqnarray*}
h_{n}(u) &:&=\wp _{n}(h(u))=\left\{ 
\begin{array}{cc}
\frac{h(u)^{{\large n}}-1}{h(u)-1} & h(u)\neq 1, \\ 
n & h(u)=1{\Large ,}%
\end{array}%
\right. \\
g_{n}(u) &:&=\wp _{n}(g(u)){\Large =}\left\{ 
\begin{array}{cc}
\frac{g(u)^{n}-1}{g(u)-1} & g(u)\neq 1, \\ 
n & g(u)=1.%
\end{array}%
\right.
\end{eqnarray*}%
Note that $g_{m}(u/n)\neq 0,$ since $g(u/n)^{{\large m}}=1$ implies $%
g_{m}(u/n)=m/n.$ Replacing $n$ by $m$ and $u$ by $u/n$ yields $(\ast )\ $%
above. Hence 
\begin{equation*}
K(h_{m}(u/n)u/n)=g_{m}(u/n)K(u/n)=g_{m}(u/n)g_{n}(u/n)^{-1}K(h_{n}(u/n)u/n),
\end{equation*}%
giving $(\ast \ast )$ above.\hfill $\square $

\bigskip

By Lemma 4.1 above for any $m,n\in \mathbb{N}$ 
\begin{equation*}
K(h_{m}(u/n)u/n)=g_{m}(u/n)K(u/n).
\end{equation*}%
Taking $m=kn$ and, separately, $m=\bar{k}n$ gives%
\begin{eqnarray*}
K(uh_{kn}(u/n)/n) &=&g_{kn}(u/n)K(u/n), \\
K(uh_{\bar{k}n}(u/n)/n) &=&g_{\bar{k}n}(u/n)K(u/n).
\end{eqnarray*}%
Eliminating $K(u/n)$ gives, as $g_{\bar{k}n}(u/n)\neq 0$ (see above),%
\begin{equation*}
K(uh_{kn}(u/n)/n)=g_{kn}(u/n)g_{\bar{k}n}(u/n)^{-1}K(uh_{\bar{k}n}(u/n)/n).
\end{equation*}%
We will deduce the radiality property from this equation by studying the
sets 
\begin{equation*}
Q^{g}:=\{g_{kn}(u/n):n\in \mathbb{N},k>0\},\quad Q^{h}:=\{h_{kn}(u/n):n\in 
\mathbb{N},k>0\},
\end{equation*}%
which it emerges are \textit{dense} in $[0,\infty ).$ The terms in both sets
take the form%
\begin{equation*}
q_{n}(k):=\frac{(1+\delta _{n})^{kn}-1}{n\delta _{n}}>0,\text{ }
\end{equation*}%
where respectively%
\begin{equation*}
\delta _{n}=\delta _{n}^{g}=g(u/n)-1,\text{ or }\delta _{n}^{h}=h(u/n)-1,
\end{equation*}%
both of which tend to $0$ (by the continuity of $g$ and $h$ at $0).$ We put%
\begin{equation*}
Q:=\{q_{n}(k):n\in \mathbb{N},k>0\}.
\end{equation*}

Writing $m=kn$ and noting that for $\delta _{n}\neq 0$%
\begin{equation*}
\frac{(1+\delta _{n})^{m}-1}{n\delta _{n}}=\frac{m}{n}+\frac{\delta _{n}}{n}%
c_{2}^{m}+...+\frac{\delta _{n}^{m-1}}{n},
\end{equation*}%
we may again use the \textit{L'Hospital convention} to interpret $q_{n}(k)$
as $m/n=k$ whenever $\delta _{n}=0$. In Proposition 4.1 below, we show that $%
Q$ is dense in $[0,\infty )$ when the sequence $n\delta _{n}$ is \textit{%
convergent}. See the Remark immediately below for the significance of this
assumption in terms of differentiability. The proof of density is achieved
by identifying 
\begin{equation*}
\ell (k):=\lim_{n\rightarrow \infty }q_{n}(k),
\end{equation*}%
which emerges as a simple increasing continuous injection for $k>0.$ This
gives an immediate way of \textit{steering} $q_{n}(k)$ into approaching any
point $s$ of $\mathbb{[}0,\infty )$ by use of the inverse function $k(s)$ of 
$\ell (k)$. We call $k(s)$ a \textit{steering function}.

If $n\delta _{n}$ is divergent to $\pm \infty $, the proof for Proposition
4.2 offers a more complicated steering function for approaching through $Q$
all the points $s$ of $[0,\infty )$.

By taking limits relative to an appropriate subset $\mathbb{N}_{u}^{\prime
}\subseteq \mathbb{N}$, we may arrange that each of the sequences $\{n\delta
_{n}^{g}\}_{n\in \mathbb{N}^{\prime }}$ and $\{n\delta _{n}^{h}\}_{n\in 
\mathbb{N}^{\prime }}$ is either convergent or divergent. There is thus no
loss of generality in assuming below that $\mathbb{N}^{\prime }=\mathbb{N}$.
We apply these results in Proposition 4.3 to prove the promised radiality
property of continuous solutions $(K,g,h)$ of the $(GGE),$ and as
corollaries compute the forms that the radial link function may take.

\bigskip

\noindent \textbf{Remark.} For $g,h$ the auxiliaries of the Goldie equation
with $g(0)=h(0)=1,$ the case(s) under Proposition 4.1 below corresponds to a
differentiability assumption. Thus, for example,%
\begin{eqnarray*}
g(u/n) &=&1+\delta _{n}: \\
g^{\prime }(0) &=&\lim_{n\rightarrow \infty }\frac{g(u/n)-1}{1/n}%
=\lim_{n}n\delta _{n}=\gamma _{g}(u): \\
g(tu) &=&1+t\gamma _{g}(u)+o(t).
\end{eqnarray*}

\noindent \textbf{Proposition 4.1.} \textit{Assume }$\delta _{n}\rightarrow
0 $\textit{\ with }$\gamma :=\lim_{n}n\delta _{n}\in \mathbb{R}$\textit{.
Then}%
\begin{eqnarray*}
\ell (k) &:&=\lim_{n}q_{n}(k)=\frac{e^{k\gamma }-1}{\gamma },\text{ with 
\textit{steering} provided by } \\
k(s) &:&=\frac{\log (1+s\gamma )}{\gamma }.
\end{eqnarray*}%
\textit{The case }$\gamma =0$\textit{\ falls under the L'Hospital convention
as }$\ell (k)=k$\textit{\ with inverse }$k(s)=s$\textit{. Given these
assumptions, }$Q$\textit{\ is dense in }$\mathbb{R}_{+}.$

\bigskip

\noindent \textbf{Remarks.} The case $\gamma =0$ requires separate proof.
Since $\ell (k)$ is a (non-constant) continuous function of $k,$ its range
is an interval $J,$ and so the set $Q$ is dense in $J.$

\bigskip

\noindent \textbf{Proof of Proposition 4.1.} Put $\gamma _{n}:=n\delta _{n}.$

\textit{Case 1.} $\gamma _{n}\rightarrow \gamma \neq 0.$ Here%
\begin{equation*}
q_{n}(k):=\frac{(1+\frac{\gamma _{n}}{n})^{kn}-1}{\gamma _{n}}\rightarrow
\ell (k):=\frac{e^{k\gamma }-1}{\gamma }=k+\frac{1}{2}\gamma k^{2}+....
\end{equation*}%
So $\ell (0+)=0$ and $\ell ^{\prime }(k)=e^{k\gamma }>0.$ Hence $\{\ell
(k):k\geq 0\}=\mathbb{[}0,\infty ),$ so that $\{q_{n}(k):n,k>0\}$ forms a
dense set.

\bigskip

\textit{Case 2.} $\gamma _{n}=n\delta _{n}\rightarrow 0.$ Expanding $\log
(1+t)$ around $t=0$ gives 
\begin{equation*}
\log (1+t)=t\frac{1}{1+d(t)},
\end{equation*}%
for some $d(t)$ between $0$ and $t.$ Take $t=\delta _{n},$ so that $d(\delta
_{n})\rightarrow 0,$ and put%
\begin{equation*}
k_{n}:=\frac{k}{1+d(\delta _{n})}\rightarrow k,\text{ and }t_{n}:=kn\log
(1+\delta _{n}).
\end{equation*}%
Then%
\begin{equation*}
t_{n}=n\delta _{n}k_{n}.
\end{equation*}%
Likewise expanding $\exp (t)$ around $t=0:$ 
\begin{equation*}
\exp (t)=1+te^{p(t)},
\end{equation*}%
for some $p(t)$ between $0$ and $t,$ gives%
\begin{eqnarray*}
q_{n}(k) &=&\frac{(1+\delta _{n})^{kn}-1}{n\delta _{n}}=\frac{\exp (kn\log
(1+\delta _{n}))-1}{n\delta _{n}} \\
&=&\frac{t_{n}e^{p(t_{n})}}{n\delta _{n}}=k_{n}e^{p(t_{n})}\rightarrow
ke^{0}=k,
\end{eqnarray*}%
since $n\delta _{n}k_{n}\rightarrow 0.$ That is,%
\begin{equation*}
\ell (k)=k,
\end{equation*}%
as asserted, and again $Q$ is dense in $[0,\infty )$.\hfill $\square $

\bigskip

\noindent \textbf{Proposition 4.2.} \textit{Assume }$\delta _{n}\rightarrow
0 $\textit{\ with }$\lim_{n}n\delta _{n}=\pm \infty $\textit{. Then }$Q$%
\textit{\ is dense in }$[0,\infty )$\textit{\ and the steering functions}%
\begin{equation*}
k_{n}(s)=\frac{\log (1+sn\delta _{n})}{n\log (1+\delta _{n})}
\end{equation*}%
\textit{secure a sequence in }$Q$\textit{\ approaching }$s.$

\bigskip

\noindent \textbf{Proof. }In view of the denominator oddness in $\delta
_{n}, $ it suffices to consider the case when $n\delta _{n}\rightarrow
+\infty .$ Furthermore, it suffices to consider the density of%
\begin{equation*}
\bar{q}_{n}(k)=\frac{(1+\delta _{n})^{kn}}{n\delta _{n}},
\end{equation*}%
since under the current assumptions%
\begin{equation*}
\bar{q}_{n}(k)-q_{n}(k)=\frac{1}{n\delta _{n}}\rightarrow 0.
\end{equation*}%
We put 
\begin{equation*}
r_{k}(n)=\log \bar{q}_{n}(k)=kn\log (1+\delta _{n})-\log (n\delta _{n}),
\end{equation*}%
which reduces the density consideration of $Q$ in $[0,\infty )$ to that of $%
\bar{Q}:=\{r_{k}(n):n\in \mathbb{N},k>0\}$ in $\mathbb{R}.$ For any $r\in 
\mathbb{R}$, consider any $n$ with $r+\log n\delta _{n}>0.$ Solving $%
r=r_{n}(k)$ for $k$ gives a steering function 
\begin{equation*}
k=\kappa _{n}(r):=\frac{(r+\log n\delta _{n})}{n\log (1+\delta _{n})}>0.
\end{equation*}%
That is, any $r\in \mathbb{R}$ is \textit{achieved} as being exactly $%
r_{k}(n)\in \bar{Q}$ for any $n$ by the choice $k=\kappa _{n}(r)>0,$ since%
\begin{equation*}
r_{k}(n)=\log \bar{q}_{n}(\kappa _{n}(r))=\frac{(r+\log n\delta _{n})}{n\log
(1+\delta _{n})}n\log (1+\delta _{n})-\log (n\delta _{n})=r.
\end{equation*}%
Furthermore, as $\kappa _{n}(r)$ is increasing in $r,$ setting $r^{\prime
}:=r+1/n,$ with $n$ large enough so that $r^{\prime }+\log n\delta _{n}>0,$
gives $\kappa _{n}(r^{\prime })>0$ and 
\begin{eqnarray*}
r_{n}(\kappa _{n}(r^{\prime })-r_{n}(\kappa _{n}(r)) &=&(\kappa
_{n}(r^{\prime })-\kappa _{n}(r))n\log (1+\delta _{n}) \\
&=&\frac{(r^{\prime }-r)}{n\log (1+\delta _{n})}n\log (1+\delta _{n})=\frac{1%
}{n}\rightarrow 0.
\end{eqnarray*}%
It now follows that any $r\in \mathbb{R}$ is achieved as a limit of points $%
r_{k}(n)$ by taking $k=\kappa _{n}(r+1/n)$ with $n$ large enough. Hence $%
\bar{Q}$ is dense in $\mathbb{R}$, and so $Q$ is dense in $\mathbb{R}_{+}.$
In particular, for $s>0$ and%
\begin{equation*}
r=\log (s+\frac{1}{n\delta _{n}}),
\end{equation*}%
we have 
\begin{equation*}
\log (s+\frac{1}{n\delta _{n}})=\log \bar{q}_{n}(\kappa _{n}(r)):\qquad s=%
\bar{q}_{n}(\kappa _{n}(r))-\frac{1}{n\delta _{n}}=q_{n}(\kappa _{n}(r)).
\end{equation*}%
Putting%
\begin{equation*}
k_{n}(s):=\kappa _{n}(\log (s+\frac{1}{n\delta _{n}})=\frac{\log (s+\frac{1}{%
n\delta _{n}})+\log n\delta _{n}}{n\log (1+\delta _{n})}=\frac{\log
(sn\delta _{n}+1)}{n\log (1+\delta _{n})}
\end{equation*}%
provides steering of the sequence $q_{n}(k_{n}(s))$ in $Q$ to the limit $s.$%
\hfill $\square $

\bigskip

\noindent \textbf{Remark.} Above, $k_{n}(r)\rightarrow 0$ as $n\rightarrow
\infty .$ To see this, note that $n\log (1+\delta _{n})\simeq n\delta _{n}$
to first order (cf. Case 2 in Prop. 4.1).

\bigskip

Proposition 4.3 is the next step in identifying the radiality property
explicitly, by reference to a function which we temporarily denote by $\bar{%
\lambda}_{u}(s)$ and which we subsequently show in Cor. 4.1 below to be the
link function $\lambda _{u}(s)$ of \S 3.

\bigskip

\noindent \textbf{Proposition 4.3.} \textit{For continuous }$(K,h,g)$\textit{%
\ solving }$(GGE),$\textit{\ such that }$n\delta _{n}^{g}$ \textit{and }$%
n\delta _{n}^{h}$\textit{\ are each either convergent or divergent
sequences, if }$K(u)\neq 0,\mathit{\ }$\textit{then there exists }$\bar{%
\lambda}_{u}(s)\geq 0$\textit{\ defined for }$s>0$\textit{\ with}%
\begin{equation*}
K(su)=\bar{\lambda}_{u}(s)K(u).
\end{equation*}

\noindent \textbf{Proof.} Fix $u$ with $K(u)\neq 0.$ By Lemma 4.1,%
\begin{equation*}
K(h_{m}(u/n)u/n)=g_{m}(u/n)K(u/n).
\end{equation*}%
Recall that%
\begin{equation*}
\delta _{n}=g(u/n)-1=\delta _{n}^{g}\text{ resp. }h(u/n)-1=\delta _{n}^{h}
\end{equation*}%
gives rise to%
\begin{equation*}
q_{n}^{g}(k)=g_{kn}(u/n)/n,\text{ resp. }q_{n}^{h}(k)=h_{kn}(u/n)/n.
\end{equation*}%
Taking $m=kn$ and separately $m=\bar{k}n$ gives%
\begin{eqnarray*}
K(uh_{kn}(u/n)/n) &=&g_{kn}(u/n)K(u/n), \\
K(uh_{\bar{k}n}(u/n)/n) &=&g_{\bar{k}n}(u/n)K(u/n),
\end{eqnarray*}%
and so%
\begin{equation*}
K(uh_{kn}(u/n)/n)=g_{kn}(u/n)K(u/n)=g_{kn}(u/n)g_{\bar{k}n}(u/n)^{-1}K(uh_{%
\bar{k}n}(u/n)/n).
\end{equation*}%
Applying one or other of Props 4.1 and 4.2 and taking the relevant steering
functions, write $k=k(s)$ (resp. $k_{n}(s)$) with $h_{kn}(u/n)/n\rightarrow
s>0,$ and likewise write $\bar{k}$ (resp. $\bar{k}_{n}(1)$) with $h_{\bar{k}%
n}(u/n)/n\rightarrow 1.\ $By continuity of $K,$%
\begin{equation*}
K(su)=\lim_{n}K(uh_{kn}(u/n)/n)\text{ and }K(u)=\lim_{n}K(uh_{\bar{k}%
n}(u/n)/n).
\end{equation*}%
So 
\begin{equation*}
K(su)=\bar{\lambda}_{u}(s)K(u),
\end{equation*}%
where the limit%
\begin{equation*}
\bar{\lambda}_{u}(s):=\lim_{n}g_{kn}(u/n)g_{\bar{k}n}(u/n)^{-1}=%
\lim_{n}q_{n}^{g}(k)q_{n}^{g}(\bar{k})^{-1}
\end{equation*}%
exists, as $K(u)\neq 0.$ Indeed, $q_{n}^{g}(k)q_{n}^{g}(\bar{k})^{-1}$
remain bounded over $n$, otherwise $K(su)$ is undefined. So $\lim_{n\in 
\mathbb{M}}q_{n}^{g}(k)q_{n}^{g}(\bar{k})^{-1}$ is identical for each
infinite $\mathbb{M\subseteq N}$.

Note that 
\begin{equation*}
q_{n}^{g}(k)q_{n}^{g}(\bar{k})^{-1}=\frac{g(u/n)^{{\large kn}}-1}{g(u/n)^{%
\bar{k}n}-1}>0\text{ or }=\frac{k}{\bar{k}}>0,
\end{equation*}%
the latter in case $g(u/n)=1.$\hfill $\square $

\bigskip

\noindent \textbf{Conclusions from density considerations}

Our next result recovers from Proposition 4.1 (i.e. the two cases $\lim n{%
\delta }_{n}$ zero or not) all the defining clauses of the link function of
Lemma 3.1, giving explicit form to the radiality result in Proposition 4.3.

\bigskip

\noindent \textbf{Corollary 4.1.} \textit{For continuous }$(K,h,g)$\textit{\
solving }$(GGE)$\textit{\ and }$u$ \textit{with }$K(u)\neq 0,$ \textit{if\ }$%
n\delta _{n}^{h}\rightarrow \gamma _{h}=\gamma _{h}(u)\in \mathbb{R}$ 
\textit{and }$n\delta _{n}^{g}\rightarrow \gamma _{g}=\gamma _{g}(u)\in 
\mathbb{R},$\textit{\ then }%
\begin{equation*}
\bar{\lambda}_{u}(s)=\left\{ 
\begin{array}{cc}
\lbrack (1+s\gamma _{h})^{\gamma _{g}/\gamma _{h}}-1]/[(1+\gamma
_{h})^{\gamma _{g}/\gamma _{h}}-1], & \text{if }\gamma _{g}\neq 0\text{ and }%
\gamma _{h}\neq 0, \\ 
(e^{s\gamma _{g}}-1)/(e^{\gamma _{g}}-1), & \text{if }\gamma _{g}\neq 0\text{
and }\gamma _{h}=0,\text{ } \\ 
\log (1+s\gamma _{h})/\log (1+\gamma _{h}), & \text{if }\gamma _{h}\neq 0%
\text{ and }\gamma _{g}=0, \\ 
s, & \text{if }\gamma _{h}=0=\gamma _{g},%
\end{array}%
\right.
\end{equation*}%
\textit{so that generally if }$\gamma _{h}=\gamma _{g},$ \textit{then}%
\begin{equation*}
\bar{\lambda}_{u}(s)=s.
\end{equation*}%
\textit{\ In all these cases }$\bar{\lambda}_{u}(s)=\lambda (s;\gamma
_{h}(u),\gamma _{g}(u))$ \textit{and is differentiable.}

\bigskip

\noindent \textbf{Proof.} We argue by cases. First suppose $\gamma _{g}\neq
0 $ and $\gamma _{h}\neq 0.$ As above $(1+\delta _{n}^{h})^{n}\rightarrow
e^{\gamma _{h}}$ and $(1+\delta _{n}^{g})^{n}\rightarrow e^{\gamma _{g}}$
and when as in Prop. 4.3 $k,\bar{k}$ are constants, we use logarithmic
steering:%
\begin{equation*}
\text{ }k^{h}=k^{h}(s)=\frac{1}{\gamma _{h}}\log (1+s\gamma _{h})\text{ and }%
\bar{k}^{h}=k^{h}(1)=\frac{1}{\gamma _{h}}\log (1+\gamma _{h}).
\end{equation*}%
So, as $\gamma _{g}\neq 0,$%
\begin{equation*}
\lambda _{u}(s)=\frac{q_{n}^{g}(k^{h})}{q_{n}^{g}(\bar{k}^{h})}=\frac{%
e^{k^{h}(s)\gamma _{g}}-1}{e^{\bar{k}^{h}(1)\gamma _{g}}-1}=\frac{(1+s\gamma
_{h})^{\gamma _{g}/\gamma _{h}}-1}{(1+\gamma _{h})^{\gamma _{g}/\gamma
_{h}}-1}.
\end{equation*}%
as required. Likewise, if $\gamma _{g}\neq 0$ and $\gamma _{h}=0,$ then%
\begin{equation*}
\lambda _{u}(s)=\frac{e^{k^{h}(s)\gamma _{g}}-1}{e^{\bar{k}^{h}(1)\gamma
_{g}}-1}=\frac{e^{s\gamma _{g}}-1}{e^{\gamma _{g}}-1},
\end{equation*}%
giving the second case.

As in Prop. 4.1 Case 2,%
\begin{eqnarray*}
q_{n}(k) &=&\frac{(1+\delta _{n})^{kn}-1}{n\delta _{n}}=\frac{\exp (kn\log
(1+\delta _{n}))-1}{n\delta _{n}} \\
&=&\frac{t_{n}e^{p(t_{n})}}{n\delta _{n}}=k_{n}e^{p(t_{n})}\rightarrow
ke^{0}=k.
\end{eqnarray*}%
So suppose now that $\gamma _{h}\neq 0$ but $\gamma _{g}=0.$ To compute $%
\lambda $ we need a combination of $k^{h}$ with $\delta _{n}^{g}.$ Here with 
$t_{n}^{g}:=k^{h}n\log (1+\delta _{n}^{g})$%
\begin{eqnarray*}
q_{n}^{g}(k^{g}) &=&\frac{(1+\delta _{n}^{g})^{k^{g}n}-1}{n\delta _{n}^{g}}=%
\frac{\exp (k^{h}n\log (1+\delta _{n}^{g}))-1}{n\delta _{n}^{g}} \\
&=&\frac{t_{n}^{g}e^{p(t_{n}^{g})}}{n\delta _{n}^{g}}=k^{g}e^{p(t_{n}^{g})}%
\rightarrow k^{g}e^{0}=k^{g}.
\end{eqnarray*}%
Similarly, with $t_{n}^{g}=n\delta _{n}^{g}k_{n}^{h}$%
\begin{equation*}
\exp (k^{h}n\log (1+\delta _{n}^{g}))-1=n\delta
_{n}^{g}k_{n}^{h}e^{p(t_{n}^{g})}\text{ },
\end{equation*}%
\begin{equation*}
k_{n}^{h}:=k^{h}/(1+d(\delta _{n}))=\frac{\log (1+s\gamma )/\gamma }{%
1+d(\delta _{n})}\rightarrow k^{h}=\log (1+s\gamma _{h})/\gamma _{h}.
\end{equation*}%
From here,%
\begin{equation*}
\lambda _{u}(s):=\lim_{n}\frac{(1+\delta _{n}^{g})^{kn}-1}{(1+\delta
_{n}^{g})^{\bar{k}n}-1}=\frac{n\delta _{n}^{g}k_{n}^{h}e^{p(t_{n}^{g})}}{%
n\delta _{n}^{g}\bar{k}_{n}^{h}e^{p(\bar{t}_{n}^{g})}}=\frac{%
k^{h}(s)e^{p(t_{n}^{g})}}{k^{h}(1)e^{p(\bar{t}_{n}^{g})}}\rightarrow \frac{%
\log (1+s\gamma _{h})}{\log (1+\gamma _{h})},
\end{equation*}%
as $t_{n}^{g}\rightarrow 0,$ since $n\delta _{n}^{g}\rightarrow \gamma
_{g}=0 $ here.

Finally, suppose $\gamma _{h}=\gamma _{g}=0.$ Then $k^{h}(s)=s$ and we have%
\begin{equation*}
\lambda _{u}(s)=s,
\end{equation*}%
as required in the last case$.$\hfill $\square $

\bigskip

We now consider the cases in which one of the sequences is divergent. For
these, it is helpful in the context of $(GGE),$ and hence also of $(GFE),$
to note that differentiability at $0$ of $K\ $in direction $u$ for $h(u)\neq
0$ implies differentiability of $K$ elsewhere along $u$, as the following
shows.%
\begin{eqnarray*}
K(u+h(u)tu)-K(u) &=&g(u)K(tu), \\
\frac{K(u+h(u)tu)-K(u)}{th(u)} &=&\frac{g(u)}{h(u)}\frac{K(tu)-K(0)}{t}.
\end{eqnarray*}%
Given this observation, the connection between $g$ and $h$ in $(GGE)\ $is
such that either one of them is differentiable at the origin iff the other
is. Indeed, starting from%
\begin{equation*}
K([t+h(tu)]u)=K(tu+h(tu)u)=K(tu)+g(tu)K(u),
\end{equation*}%
it follows, after subtracting $K(u)$ from each side (and since $K(0)=0),$
that%
\begin{eqnarray*}
&&\frac{K((t+h(tu))u)-K(h(0)u)}{(t+h(tu)-1)}\cdot \frac{t+(h(tu)-1)}{t} \\
&=&\frac{K(tu)-K(0)}{t}+\frac{(g(tu)-1)}{t}K(u),
\end{eqnarray*}%
for all small enough $t>0.$ Given radiality and so differentiability of $%
K_{u}$ established in Cor. 4.1, the preceding equation also yields, for
non-zero $K(u),$ 
\begin{equation*}
\lambda _{u}^{\prime }(1)\cdot (1+\lim_{t\rightarrow 0}\frac{(h(tu)-1)}{t}%
)=\lambda _{u}^{\prime }(0)+\lim_{t\rightarrow 0}\frac{(g(tu)-1)}{t},
\end{equation*}%
implying, since $\lambda _{u}^{\prime }>0$, that $n\delta
_{n}^{h}\rightarrow \infty $ iff $n\delta _{n}^{g}\rightarrow \infty ,$ as
is borne out below.

\bigskip

\noindent \textbf{Corollary 4.2}. \textit{For continuous }$(K,h,g)$\textit{\
solving }$(GGE)$ \textit{and }$u$ \textit{with }$K(u)\neq 0,$ 
\begin{equation*}
n\delta _{n}^{g}\rightarrow \infty \text{ iff }n\delta _{n}^{h}\rightarrow
\infty ,
\end{equation*}%
\textit{and for }$s>0$%
\begin{equation*}
\lambda _{u}(s)=s.
\end{equation*}%
\textit{Thus }$\lambda _{u}$ \textit{is differentiable in all cases and so }$%
K_{u}$ \textit{is differentiable.}

\bigskip

\noindent \textbf{Proof.} First we suppose that $n\delta _{n}^{h}\rightarrow
\infty .$ Here, according to Prop. 4.2, steering through $Q$ towards $s>0$
is provided by%
\begin{equation*}
k_{n}^{h}(s)=\frac{\log (1+sn\delta _{n}^{h})}{n\log (1+\delta _{n}^{h})}.
\end{equation*}%
So%
\begin{eqnarray*}
\log \bar{q}_{n}(k_{n}^{h}(s)) &=&\frac{\log (s+\frac{1}{n\delta _{n}^{h}}%
)+\log n\delta _{n}^{h}}{n\log (1+\delta _{n}^{h})}n\log (1+\delta
_{n}^{h})-\log (n\delta _{n}^{h})=\log (s+\frac{1}{n\delta _{n}^{h}}): \\
\bar{q}_{n}(k_{n}^{h}(s)) &=&s+\frac{1}{n\delta _{n}^{h}}.
\end{eqnarray*}%
We now consider that%
\begin{equation*}
\frac{(1+\delta _{n}^{g})^{kn}-1}{n\delta _{n}^{g}}=q_{n}^{g}(k_{n}^{h})=%
\bar{q}_{n}^{g}(k_{n}^{h})-\frac{1}{n\delta _{n}^{g}}=(s+\frac{1}{n\delta
_{n}^{h}})-\frac{1}{n\delta _{n}^{g}}.
\end{equation*}%
Hence, with $k$ for $k^{h}(s)$ and $\bar{k}$ for $k^{h}(1),$ 
\begin{eqnarray*}
\lambda (s) &=&\lim_{n}\frac{q_{n}^{g}(k)}{q_{n}^{g}(\bar{k})}=\lim_{n}\frac{%
(1+\delta _{n}^{g})^{{\large kn}}-1}{(1+\delta _{n}^{g})^{{\large \bar{k}n}%
}-1}=\lim_{n}\frac{n\delta _{n}^{g}(s+\frac{1}{n\delta _{n}^{h}})-1}{n\delta
_{n}^{g}(1+\frac{1}{n\delta _{n}^{h}})-1} \\
&=&\lim_{n}\frac{(s+\frac{1}{n\delta _{n}^{h}})-1/(n\delta _{n}^{g})}{(1+%
\frac{1}{n\delta _{n}^{h}})-1/(n\delta _{n}^{g})}.
\end{eqnarray*}%
So, provided $\gamma _{g}\neq 1,$ we conclude that since $n\delta
_{n}^{h}\rightarrow +\infty ,$%
\begin{equation*}
\lambda _{u}(s)=\frac{\gamma _{g}s-1}{\gamma _{g}-1}\text{ or }s,\text{
according as }n\delta _{n}^{g}\rightarrow \gamma _{g}\in \mathbb{R}\text{ or 
}n\delta _{n}^{g}\rightarrow \infty .\text{ }
\end{equation*}%
Taking $s\neq 1,$ Proposition 4.3 implies that the case $n\delta
_{n}^{g}\rightarrow \gamma _{g}=1$ cannot arise as $K(u)\neq 0$ yields the
finiteness of $\lambda _{u}(s)>0$. We `park this case' temporarily, but will
eventually show that also $n\delta _{n}^{g}\rightarrow \gamma _{g}\in 
\mathbb{R}$ cannot occur.

Now suppose that $n\delta _{n}^{h}\rightarrow \gamma _{h}\in \mathbb{R}$ but 
$n\delta _{n}^{g}\rightarrow +\infty .$ Taking%
\begin{equation*}
k^{h}(s)=\frac{\log (1+s\gamma _{h})}{\gamma _{h}}
\end{equation*}%
gives%
\begin{equation*}
\lambda _{u}(s)=\lim_{n}\frac{q_{n}^{g}(k^{h})}{q_{n}^{g}(\bar{k}^{h})}%
=\lim_{n}\frac{\bar{q}_{n}^{g}(k^{h})}{\bar{q}_{n}^{g}(\bar{k}^{h})}=\lim_{n}%
\frac{(1+\delta _{n}^{g})^{kn}}{(1+\delta _{n}^{g})^{\bar{k}n}}.
\end{equation*}%
So again, since $\lambda _{u}(s)>0$ as $K(u)\neq 0,$%
\begin{equation*}
\log \lambda _{u}(s)=\left( \frac{\log (1+s\gamma ^{b})/(1+\gamma ^{b})}{%
\gamma ^{b}}\right) n\log (1+\delta _{n}^{g}).
\end{equation*}%
Here the RHS\ is divergent, since as above%
\begin{equation*}
n\log (1+\delta _{n}^{g})=\frac{n\delta _{n}^{g}}{1+d(\delta _{n}^{g})}%
\rightarrow \infty .
\end{equation*}%
Again, by Prop. 4.3, this case cannot arise.

But now that all forms of $\lambda _{u}(s)$ are known, we see that $\lambda
_{u}$ is differentiable and, by the radiality property, so is $K_{u}.$

It now follows that in the `parked case' with $n\delta _{n}^{h}\rightarrow
+\infty $ and $n\delta _{n}^{g}\rightarrow \gamma _{g}\in \mathbb{R}$, in
fact $g$ is differentiable at $0.$ So by our initial observations in the
introductory paragraph, $h$ is differentiable at $0,$ contrary to $n\delta
_{n}^{h}\rightarrow +\infty .$ That is, the first case stated in the parked
case does not in fact arise.\hfill $\square $

\bigskip

\noindent \textbf{Remark}. The case $n\delta _{n}^{h}\rightarrow \gamma _{h}$
but $n\delta _{n}^{g}\rightarrow +\infty $ above, ruled out by Prop. 4.3, is
also ruled out by the fact that $h$ and $K$ being differentiable implies $g$
is differentiable.

\bigskip

Our next result addresses the partially pexiderized $(GS)\ $functional
equation below, whose solution we will need later.%
\begin{equation*}
g(su+h(su)tu)=g(su)g(tu).
\end{equation*}%
For $g>0,$ taking $\kappa (tu)=\log g(tu)$ gives rise to a special case of $%
(GGE)$ identified below (with $\kappa $ replacing $K$).

\bigskip

\noindent \textbf{Proposition 4.4.} \textit{For continuous }$(K,h,g)$\textit{%
\ solving }$(GGE),$ \textit{the kernel }$K$ \textit{is differentiable along }%
$u$\textit{\ for }$h(u)>0,$\textit{\ }$K(u)\neq 0$\textit{\ and is strictly
increasing along }$u$\textit{. Hence if }$K(u)=0,$ \textit{and }$h(su)>0$ 
\textit{with }$s>0,$ \textit{then }%
\begin{equation*}
K(su)=0,\text{ i.e. }K(su)=\lambda _{u}(s)K(u)\text{ with }\lambda
_{u}(s)\equiv 0.
\end{equation*}%
\textit{Furthermore, if either auxiliary generates a Popa binary operation
along }$u,$ \textit{then so does the other.}

\textit{In particular, continuous solutions of the equation}%
\begin{equation*}
\kappa (a+h(a)b)=\kappa (a)+\kappa (b)\text{ for }a,b\in \langle u\rangle
\end{equation*}%
\textit{have }$h(tu)\equiv 1+t\rho u$\textit{\ for some }$\rho \in \mathbb{R}
$\textit{.}

\bigskip

\noindent \textbf{Proof.} The first statement follows from Corollaries 4.1
and 4.2 and Prop 4.3. Suppose that $K(u)=0.$ We claim that $K(su)=0$ for all 
$s>0.$ Suppose not and that $K(ru)\neq 0$ for some $r>0.$ Then with $%
s=r^{-1} $%
\begin{equation*}
0=K(u)=K(sru)=\lambda _{ru}(r^{-1})K(ru)\neq 0.
\end{equation*}%
Since by Corollaries 4.1 and 4.2, $\lambda _{ru}(r^{-1})>0,$ this is a
contradiction. So $K(ru)=0$ for all $r>0.$

Putting $U=K(u),V=K(v),H=hK^{-1}$ and $G=gK^{-1},$%
\begin{eqnarray*}
K(u+h(u)v) &=&K(u)+g(u)K(v)]: \\
U+H(U)V &=&K^{-1}[U+G(U)V].
\end{eqnarray*}%
Suppose $h(u)=1+\rho u.$ Then, as $u+h(u)v=u+v+\rho uv$ is commutative in $%
u,v,$%
\begin{eqnarray*}
K(u\circ _{\rho }v) &=&K(u)+g(u)K(v)=K(v)+g(v)K(u), \\
(g(u)-1)K(v) &=&K(u)(g(v)-1), \\
K(v)/(g(v)-1) &=&K(u)/(g(u)-1)=c,\text{ say,} \\
g(u) &=&1+cK(u), \\
G(x) &=&g(K^{-1}(x))=1+cx.
\end{eqnarray*}%
Conversely, suppose $G(U)=1+cU.$ Then as above $U+G(U)V=U+V+cUV,$ so%
\begin{eqnarray*}
U+H(U)V &=&K^{-1}[U+G(U)V]=V+H(V)U, \\
(H(U)-1)V &=&(H(V)-1)U, \\
(H(U)-1)/U &=&(H(V)-1)/V=\rho ,\text{ say,} \\
H(U) &=&1+\rho U.
\end{eqnarray*}%
In either case commutativity implies that the auxiliary on the other side of
the equation generates a binary group operation.

The final assertion arises by a specialization of $(GGE)\ $to $g(x)\equiv 1,$
permitted by Prop. 4.3.\hfill $\square $

\bigskip

\noindent \textbf{Corollary 4.3.} (i) \textit{For }$s,t>0,$%
\begin{equation*}
\lambda _{su}(t)=\lambda _{u}(st)/\lambda _{u}(s)\text{ and }\lambda
_{su}^{\prime }(t)=\lambda _{u}^{\prime }(st)s/\lambda _{u}(s).
\end{equation*}%
(ii) \textit{For }$u$\textit{\ with }$h(u)>0$\textit{\ and }$h(-u)>0:$%
\begin{equation*}
K(-u)=-g(-u)\lambda _{u}(1/h(-u))K(u).
\end{equation*}

\noindent \textbf{Proof.} (i) As for the first assertion, $\lambda
_{u}(st)=\lambda _{su}(t)\lambda _{u}(s),$ this holds, since%
\begin{equation*}
K(tsu)=\lambda _{u}(st)K(u)=\lambda _{su}(t)K(su)=\lambda _{su}(t)\lambda
_{u}(s)K(u),\text{ for }s,t\geq 0.
\end{equation*}%
The second assertion follows from Prop. 4.4 by differentiation w.r.t. $t.$

(ii) This is immediate from%
\begin{equation*}
0=K(-u+h(-u)u/h(-u))=K(-u)+g(-u)K(u/h(-u)).
\end{equation*}%
\hfill $\square $

\bigskip

\noindent \textbf{5. Proof of Theorem 3.1. }

Fix $u$ and $t.$ Recall from Prop. 4.3 that%
\begin{equation*}
K(uh_{kn}(u/n)/n)=g_{kn}(u/n)g_{\bar{k}n}(u/n)^{-1}K(uh_{\bar{k}n}(u/n)/n).
\end{equation*}%
Again as in Prop. 4.3 with $k=t,$ $\bar{k}=1$ and integers $m(n)$ with $%
m(n)/n\rightarrow t$ writing $\gamma (u)$ for $\gamma _{g}(u),$ as in Prop
4.1 and Prop 4.2,%
\begin{equation*}
\frac{g(u/n)^{tn}-1}{g(u/n)-1}\rightarrow \frac{e^{t\gamma (u)}-1}{\gamma (u)%
}\text{and }\frac{g(u/n)^{m(n)}-1}{g(u/n)^{n}-1}\rightarrow \frac{e^{t\gamma
(u)}-1}{e^{\gamma (u)}-1},
\end{equation*}%
appropriately interpreted for $\gamma =0.$ Similarly, for $h$ and writing $%
\rho (u)$ for $\gamma _{h}(u),$%
\begin{equation*}
\frac{h(u/n)^{tn}-1}{h(u/n)-1}\rightarrow \frac{e^{t\rho (u)}-1}{\rho (u)}.
\end{equation*}%
By continuity of $K,$ these limit relations lead to the first claims of
Theorem 3.1, that%
\begin{equation*}
K\left( \frac{e^{t\rho (u)}-1}{\rho (u)}u\right) =\frac{e^{t\gamma (u)}-1}{%
e^{\gamma (u)}-1}K\left( \frac{e^{t\gamma (u)}-1}{\gamma (u)}\right) .
\end{equation*}%
To proceed further, reparametrize the coefficient of $u$ by taking $%
s:=(e^{t\rho (u)}-1)/\rho (u)$ and let $s=\mu =\mu (u):=(e^{\rho
(u)}-1)/\rho (u)$ correspond to $t=1.$ Solving for $t$ in terms of $s$ gives 
$e^{t\rho (u)}=(1+s\rho (u))$ and so%
\begin{equation*}
K(su)=\frac{(1+\rho (u)s)^{\gamma (u)/\rho (u)}-1}{e^{\gamma (u)}-1}K(\mu
u):\quad K(u)=\frac{(1+\rho (u))^{\gamma (u)/\rho (u)}-1}{e^{\gamma (u)}-1}%
K(\mu u),
\end{equation*}%
the latter by specializing the former to $s=1$. After cross-substitution,%
\begin{equation*}
K(su)=\frac{(1+s\rho (u))^{\gamma (u)/\rho (u)}-1}{(1+\rho (u))^{\gamma
(u)/\rho (u)}-1}K(u)=\lambda _{u}(s,\rho (u),\gamma (u))K(u).
\end{equation*}%
Here if $\gamma (u)=\rho (u),$ then $K(su)=sK(u)$ and then $%
\{u:K(su)=sK(u)\}=\{u:\gamma (u)=\rho (u)\}$ is a linear subspace.\hfill $%
\square $

\bigskip

\noindent \textbf{6. Shuffling and switching}

The defining formula $(\dag )$ (\S 3, Definition 2) of the link function
involves all the standard homomorphisms between different \textit{scalar}
Popa groups, i.e. Popa groups on $\mathbb{R}$ with $(\eta _{\rho })$ in \S 1
specialized to $\rho (t)=rt$ for $r,t\in \mathbb{R}$, so between $\mathbb{G}%
_{r}$ and $\mathbb{G}_{\theta }$ for\textit{\ }$r,\theta \in \lbrack
0,\infty ]$. These are summarized in the table below for $\mathbb{G}_{r}=%
\mathbb{G}_{r}(\mathbb{R}_{+})$ and $\mathbb{G}_{\theta }=\mathbb{G}_{\theta
}(\mathbb{R}),$ reproduced for convenience from [BinO8] Theorem BO (cf.
[Ost3]; [BinO5] [BinO6]).

\renewcommand{\arraystretch}{1.25}%
\begin{equation*}
\begin{tabular}{|l|l|l|l|}
\hline
Popa parameter & $\theta =0$ & $\theta \in (0,\infty )$ & $\theta =\infty $
\\ \hline
$r=0$ & $\kappa t$ & $\eta _{\theta }^{-1}(e^{\theta \kappa t})$ & $%
e^{\kappa t}$ \\ \hline
$r\in (0,\infty )$ & $\log \eta _{r}(t)^{\kappa /r}$ & $\eta _{\theta
}^{-1}(\eta _{r}(t)^{\theta \kappa /r})$ & $\eta _{r}(t)^{\kappa /r}$ \\ 
\hline
$r=\infty $ & $\log t^{\kappa }$ & $\eta _{\theta }^{-1}(t^{\theta \kappa })$
& $t^{\kappa }$ \\ \hline
\end{tabular}%
\end{equation*}

\renewcommand{\arraystretch}{1}

Theorem 3.1 may now be interpreted as a \textit{shuffling}, via the link
function, of these Popa homomorphisms. Explicitly, Theorem 3.1 may be read
as saying (see Corollary 6.1 below) that for a given pair $(K,g)$ satisfying 
$(GFE)$, the kernel $K$ induces a map between the canonical scalar Popa
homomorphisms.

\bigskip

\noindent \textbf{Corollary 6.1. }\textit{For }$(K,g)$\textit{\ satisfying }$%
(GFE$-$\rho ),\gamma =\log g,$\textit{\ and fixed }$u\in X,$ \textit{with }$%
\rho (u)>0$\textit{\ and }$\gamma (u)=1$\textit{, put}%
\begin{equation*}
a(u):=e^{\rho (u)}-1.
\end{equation*}%
\textit{Then:}

\noindent (i) \textit{The map} $a:\mathbb{G}_{\rho }(X)\rightarrow \mathbb{G}%
_{I}(\mathbb{R})$\textit{\ is additive, for }$I(t):=t$\textit{\ }$(t\in 
\mathbb{R)};$\textit{\ equivalently, with }%
\begin{equation*}
b(u):=\log (1+a(u)),
\end{equation*}%
$b$ \textit{is linear on }$\mathbb{G}_{\rho }(X):$\textit{\ }%
\begin{equation*}
a(u+v)=a(u)\circ _{I}a(v),\qquad b(u+v)=b(u)+b(v).
\end{equation*}

\noindent (ii)\textit{\ Put }$c(u):=\gamma (u)/\rho (u);$ \textit{then up to
the constant factor }$\eta _{a(u)}^{-1}(\eta _{a(u)}(1)^{c(u)})$\textit{\
below,}%
\begin{equation*}
\mathit{\ }K_{u}:\langle u\rangle _{\rho }\sim \mathbb{G}_{a(u)}(\mathbb{R}%
)\rightarrow \langle K(u)\rangle _{a(u)}
\end{equation*}%
\textit{induces a map} $\mathbb{G}_{a(u)}(\mathbb{R})\rightarrow \mathbb{G}%
_{a(u)}(\mathbb{R}):$%
\begin{equation*}
\eta _{a(u)}^{-1}(\eta _{a(u)}(1)^{c(u)})\cdot K(su)=\eta _{a(u)}^{-1}(\eta
_{a(u)}(s)^{c(u)})K(u).
\end{equation*}%
\textit{So again} 
\begin{equation*}
K(\langle u\rangle _{\rho })\subseteq \langle K(u)\rangle _{Y}.
\end{equation*}

\noindent (iii)\textit{\ Taking, for any }$w\in \mathbb{G}_{\rho }(X),$%
\textit{\ }%
\begin{equation*}
b_{K}(w):=e^{c(w)\log [1+\rho (w)]}-1,\qquad \psi _{\rho (w)}(t):=\eta
_{\rho (w)}^{-1}(e^{t\log [1+\rho (w)]}),
\end{equation*}%
\textit{the kernel function }$K\ $\textit{induces a map between
homomorphisms }$\mathbb{G}_{\rho (w)}(\mathbb{R})\rightarrow \mathbb{G}%
_{b_{K}(w)}(\mathbb{R}):$%
\begin{equation*}
K(\psi _{\rho (w)}(t)w)=\left\{ 
\begin{array}{cc}
\psi _{b_{K}(w)}(t)K(w), & c(w)\neq 0, \\ 
tK(w), & c(w)=0.%
\end{array}%
\right.
\end{equation*}

\noindent \textbf{Proof. }(i) Since $\rho $ is additive, by definition of $%
a: $ 
\begin{eqnarray*}
(1+a(u))(1+a(v)) &=&e^{\rho (u)+\rho (v)}=e^{\rho (u+v)}=1+a(u+v): \\
a(u+v) &=&a(u)+a(v)+a(u)a(v)=a(u)\circ _{1}a(v): \\
b(u)+b(v) &=&\log (1+a(u))(1+a(v))=\log (1+a(u+v))=b(u+v).
\end{eqnarray*}

\noindent (ii) For $h_{u}(t)=1+t\rho (u)$ with $h_{u}^{\prime }(t)\equiv
\rho (u)),$ take $w=w(u):=(e^{\rho (u)}-1)u/\rho (u)\in \langle u\rangle
_{\rho }$ (as $\rho (w)>0).$ By homogeneity of directional derivatives,%
\begin{equation*}
\rho (w(u))=e^{\rho (u)}-1,\qquad \gamma (w(u))=(e^{\rho (u)}-1)\gamma
(u)/\rho (u)=\rho (w(u))\gamma (u)/\rho (u):
\end{equation*}
\begin{equation*}
c(u):=\gamma (u)/\rho (u)=\gamma (w(u))/\rho (w(u))=c(w(u)).
\end{equation*}

The operation $\circ _{\rho }$ on $\langle w\rangle _{\rho }$ is the same as 
$\circ _{\alpha (u)}$ on $\mathbb{R}$, since $\rho (w(u))=e^{\rho
(u)}-1=a(u):$ indeed,%
\begin{equation*}
sw\circ _{\rho }tw=sw+tw+stw\rho (w)=[s+t+sta(u)]w=[s\circ _{a(u)}t]w.
\end{equation*}%
As $\alpha :=a(u)\neq 0,$ we may write $\eta _{\alpha }(s):=1+a(u)s,$ and
put 
\begin{equation*}
t=\frac{\log \eta _{\alpha }(s)}{\log \eta _{\alpha }(1)}=\frac{\log
[1+s(e^{\rho (u)}-1)]}{\rho (u)}:\qquad s=\frac{e^{\rho (u)t}-1}{e^{\rho
(u)}-1}.
\end{equation*}%
As $a(u):=e^{\rho (u)}-1$ and $c(u):=\gamma (u)/\rho (u),$ as in Theorem
3.1, with $w$ for $w(u):$%
\begin{eqnarray*}
K(sw(u)) &=&\frac{e^{\gamma (u)t}-1}{e^{\gamma (u)}-1}K(w)=\frac{\eta
_{\alpha }(s)^{c(u)}-1}{\eta _{\alpha }(1)^{c(u)}-1}K(w)e^{\rho (u)}-1 \\
&=&\frac{[\eta _{\alpha }(s)^{c(u)}-1]/\alpha }{[\eta _{\alpha
}(1)^{c(u)}-1]/\alpha }K(w)=\frac{\eta _{\alpha }^{-1}(\eta _{\alpha
}(s)^{c(u)})}{\eta _{\alpha }^{-1}(\eta _{\alpha }(1)^{c(u)})}K(w(u)).
\end{eqnarray*}%
\noindent (iii) With $w(u)$ as in (ii) above,%
\begin{eqnarray*}
\rho (w(u)) &=&e^{\rho (u)}-1\text{ }(=a(u)):\qquad \rho (u)=\log [1+\rho
(w(u))]: \\
\log [1+b_{K}(w)] &=&\rho (u):\qquad b_{K}(w):=e^{\log [1+\rho (w(u))]}-1.
\end{eqnarray*}%
Substitution into the formula $(Rad)$ of Theorem 3.1 yields the
assertion.\hfill $\square $

\bigskip

As a further corollary of Theorem 3.1, we now have the following radial
version of a familiar result (see e.g. [BinGT, Proof of Lemma 3.2.1],
[BinO3, Th. 1(ii)], [BojK, (2.2)], [AczG]), here written as result on 
\textit{switching} between $tu$ and $u$. (We will encounter a skeletal
version within the proof of Corollary 6.2 below.)

\bigskip

\noindent \textbf{Lemma 6.1}. \textit{For }$(K,g)$\textit{\ satisfying }$%
(GFE $-$\rho )$\textit{\ with }$K(u)\neq 0$ \textit{and with }$g\neq 1$%
\textit{\ on }$\langle u\rangle _{\rho }$ \textit{except at }$0:$%
\begin{equation*}
(g(tu)-1)K(u)=(g(u)-1)K(tu)\qquad (tu\in \mathbb{G}_{\rho }(X)),
\end{equation*}%
\textit{that is,}%
\begin{equation*}
(g(x)-1)K(u)=(g(u)-1)K(x)\qquad (x\in \langle u\rangle _{\rho }).
\end{equation*}

\bigskip

\noindent \textbf{Proof}. Here $u\neq 0$ (since $K(0)=0,$ by $(GFE)).$ As $%
\langle u\rangle _{\rho }$ is abelian ([BinO8, \S 3 Lemma]),%
\begin{equation*}
K(su\circ _{\rho }tu)=K(su)+g(su)K(tu)=K(tu)+g(tu)K(su).
\end{equation*}%
As $K(u)\neq 0$ and $g(su)\neq 1$ for $s\neq 0,$ Theorem 3.1 yields 
\begin{equation}
K(tu)=\lambda _{u}(t)K(u),  \tag{$R$}
\end{equation}%
whence%
\begin{equation*}
K(tu)[g(su)-1]=[g(tu)-1]K(su):\qquad \lambda _{u}(t)/[g(tu)-1]=\lambda
_{u}(s)/[g(su)-1].
\end{equation*}%
So this is constant, say $k(u).$ Hence%
\begin{equation*}
\lbrack g(tu)-1]K(u)=k(u)\lambda _{u}(t)K(u)=k(u)K(tu),
\end{equation*}%
again using $(R).$ Take $t=1;$ then%
\begin{equation*}
\lbrack g(u)-1]K(u)=k(u)K(u):\qquad \lbrack g(u)-1]=k(u).\qquad \qquad
\qquad \qquad \square
\end{equation*}

Lemma 6.2 secures the non-triviality of the radial function $g_{u}(t):=g(tu)$%
.

\bigskip

\noindent \textbf{Lemma 6.2. }\textit{For }$g$\textit{\ continuous
satisfying }$(M)$, \textit{if }$g(u)\neq 1$\textit{\ and }$\rho (u)=1,$%
\textit{\ then }$g(tu)\neq 1$\textit{\ for }$t\neq 0.$

\bigskip

\noindent \textbf{Proof. }From Lemma 2.1, $\gamma :=\log g$ satisfies $(GFE)$
in the simpler additive form $(A).$ So Theorem 3.1 here yields%
\begin{equation*}
\gamma \left( \frac{e^{t}-1}{e-1}u\right) =\frac{e^{\gamma (u)t}-1}{%
e^{\gamma (u)}-1}\gamma (u),
\end{equation*}%
as $\gamma (u)\neq 0.$ So for $t\neq 0,$ $\gamma (tu)\neq 0,$ and so $%
g(tu)\neq 1.$\hfill $\square $

\bigskip

Theorem 7.1 and Corollary 6.2 below will immediately imply our first main
result, Theorem 7.2 below, on the existence of $\sigma _{g}.$ As noted in
the Introduction, the dichotomy below concerning two null spaces arises
because two hyperplanes passing through the origin (representing the pair of
null spaces of interest) have intersection with co-dimension 1, when they
coincide, but 2 otherwise.

\bigskip

\noindent \textbf{Corollary 6.2.} \textit{Suppose }$(K,g)$ \textit{satisfies 
}$(GFE$-$\rho ),$ \textit{so that for some }$\alpha :\mathbb{G}_{\rho
}(X)\rightarrow \mathbb{R}$\textit{\ additive, }$\beta \in \mathbb{R}$%
\textit{, }$g$ \textit{is characterized in Theorem 2.1 as having the form}%
\begin{equation*}
g(x)=\bar{g}_{\alpha ,\beta }(x):=e^{\alpha (x)}(1+\rho (x))^{\beta }.
\end{equation*}%
\textit{\ Then the restriction of the kernel }$K|\mathcal{N}(\rho )$\textit{%
\ is linear on }$\mathcal{N}^{\ast }(\alpha )=\mathcal{N}(\alpha )\cap \{x:%
\bar{g}(x)=1\}$\textit{. Furthermore, either }%
\begin{equation}
\mathcal{N}(\rho )=\mathcal{N}(\alpha )  \tag{$N_{A}^{=}$}
\end{equation}%
\textit{holds, or else }$K|\mathcal{N}^{\ast }(\alpha )=0,$ \textit{and then 
}%
\begin{equation}
K(\mathcal{N}(\rho ))=\langle K(v)\rangle \text{ for some }v\in \mathcal{N}%
(\rho )\backslash \mathcal{N}(\alpha ).  \tag{$N_{B}^{=}$}
\end{equation}

\noindent \textbf{Proof.} Since $N_{A}^{=}$ holds and 
\begin{equation*}
\bar{g}(x)=e^{\alpha (x)}(1+\rho (x))^{\beta },
\end{equation*}%
$g|\mathcal{N}^{\ast }(\alpha )\equiv 1$ (as $\rho (x)=0$ here), and so
additivity of $K$ on $\mathcal{N}^{\ast }(\alpha )$ and hence (by
continuity) its linearity is immediate. If $\mathcal{N}(\rho )=\mathcal{N}%
(\alpha ),$ then $K|\mathcal{N}(\rho )$ is linear. Otherwise $\mathcal{N}%
^{\ast }(\alpha )$ is of co-dimension 1 in the subspace $\mathcal{N}(\rho )$
(see e.g. [Con, 3.5.1]). In particular, we may choose and fix $v_{2}\in 
\mathcal{N}(\rho )\backslash \mathcal{N}(\alpha ).$ Now take $v_{1}\in 
\mathcal{N}^{\ast }(\alpha )$ arbitrarily. Then as $v_{1},v_{2}\in \mathcal{N%
}(\rho )$ by commutativity and $(GFE):$%
\begin{eqnarray*}
K(v_{2})+e^{\alpha (v_{2})}K(v_{1}) &=&K(v_{1}+v_{2})=K(v_{1})+e^{\alpha
(v_{1})}K(v_{2}), \\
K(v_{1})[e^{\alpha (v_{2})}-1] &=&K(v_{2})[e^{\alpha (v_{1})}-1]=0\qquad 
\text{(as }e^{\alpha (v_{1})}=1\text{):} \\
K(v_{1}) &=&0\qquad \qquad \text{(as }e^{\alpha (v_{2})}\neq 1\text{).}
\end{eqnarray*}%
That is, $K|\mathcal{N}^{\ast }(\alpha )=0,$ and, by Theorem 3.1, $K(%
\mathcal{N}(\rho ))=\langle K(v_{2})\rangle .$ \hfill $\square $

\bigskip

\noindent \textbf{7. Inducing a Popa structure in }$Y$\textbf{\ from }$%
(GFE). $

We now turn our attention to inducing a Popa-group structure on $Y$ from a
pair $(K,g)$ satisfying $(GFE).$ Recall that $\langle \Sigma \rangle $
denotes the \textit{linear span} of $\Sigma $.

Our main result, Theorem 7.2 below, is motivated by attempting an operation
on the image $K(X)$ utilizing a solution $(K,g)$ of the ($GFE)$ via: 
\begin{equation*}
y\circ y^{\prime }=y+g(x)y^{\prime }\text{ for some }x\text{ with }y=K(x),
\end{equation*}%
which faces an obstruction, unless $K(x_{1})=K(x_{2})$ implies $%
g(x_{1})=g(x_{2}),$ i.e. $g(x_{1}-x_{2})=1.$ This is resolved in the
following

\bigskip

\noindent \textbf{Theorem 7.1.} \textit{For }$(K,g)$\textit{\ satisfying }$%
(GFE$-$\rho )$ \textit{with }$g\neq 1,$ \textit{there exists }$\sigma
:Y\rightarrow \mathbb{R}$ \textit{\ such that }$g=g^{\sigma }$\textit{\ iff
one of the} \textit{following two conditions holds: }%
\begin{equation}
\mathcal{N}(\rho )\subseteq \mathcal{N}(\gamma ),  \tag{$N_{A}$}
\end{equation}%
\textit{for} $\gamma =\log g$ \textit{together with the range condition}%
\begin{equation*}
\mathcal{R}(K)\neq K(\mathcal{N}(\rho )),
\end{equation*}%
\textit{or}%
\begin{equation}
K(\mathcal{N}(\rho ))\subseteq \langle K(u)\rangle \text{ for some }u\text{
with }g(u)\neq 1.  \tag{$N_{B}$}
\end{equation}%
\textit{Then }$\sigma $\textit{\ is uniquely determined on} $K(X).$

\bigskip

\noindent \textbf{Proof. }We first establish necessity.\textbf{\ }We suppose
the pair $(K,g)$ satisfies $(GFE)$ with $g=g^{\sigma }$ for some continuous
linear $\sigma :Y\rightarrow \mathbb{R}.$ The result follows from the
Abelian Dichotomy of [BinO8, \S 6] that either\newline
(i) $K(\mathcal{N}(\rho ))\subseteq \mathcal{N}(\sigma ).$ or \newline
(ii) $K(\mathcal{N}(\rho ))\subseteq \langle K(u)\rangle _{\sigma }$ for
some $u\in X$ with $\sigma (K(u))\neq 0$. \newline
In case (i), 
\begin{equation*}
\sigma (K(u))=0\text{ for }u\in \mathcal{N}(\rho ),
\end{equation*}%
so that, for such $u,$ $g(u)=1,$ i.e. $\gamma (K(u))=0.$ Thus $\mathcal{N}%
(\rho )\subseteq \mathcal{N}(\gamma )$ and so $(N_{A})$ holds. Furthermore,
if the range condition were to fail, then $\mathcal{R}(K)=K(\mathcal{N}(\rho
))\subseteq \mathcal{N}(\sigma ),$ implying that $\sigma (K(x))=0$ for all $%
x,$ i.e. that $g=1$ and%
\begin{equation*}
K(x\circ _{\rho }y)=K(x)+K(y).
\end{equation*}

\noindent Otherwise (ii) holds, i.e. $K(\mathcal{N}(\rho ))\subseteq \langle
K(u)\rangle _{\sigma }$ for some $u\in X$ with $\sigma (K(u))\neq 0,$ so in
particular with $g^{\sigma }(u)\neq 1,$ and a fortiori $(N_{B})$ holds.

Note that $(N_{B})$ needs no subscript on $\langle K(u)\rangle $ as $K|%
\mathcal{N}(\rho ))$ is linear, so $K(\mathcal{N}(\rho ))$ is a subspace of $%
Y.$ \hfill $\square _{\text{necessity}}$

\bigskip

The converse direction requires the construction of $\sigma $ from $g,$ so
is quite involved.\textbf{\ }Theorem 7.2 below asserts uniqueness and
sufficiency, with the latter following from Propositions 7.1A and 7.1B, our
next results. This will involve \textit{complemented} subspaces.\textbf{\ }%
We note that\textbf{\ }in the context of $Y$ a Banach space, \textit{%
algebraically} complementary spaces are \textit{topologically} complementary
[Con, Th. 13.1]. See also [LinT]. We recall our blanket assumption that $%
\rho \neq 0.$

\bigskip

\noindent \textbf{Proposition 7.1A.}\textit{\ If }$(K,g)$\textit{\ satisfies
both }$(GFE$-$\rho )$\textit{\ with }$g\neq 1$ \textit{and also }$(N_{A}),$ 
\textit{that is,}%
\begin{equation*}
\mathcal{N}(\rho )\subseteq \mathcal{N}(\gamma ),
\end{equation*}%
\textit{for }$\gamma =\log g$ \textit{(on }$\mathbb{G}^{\ast }(X)$\textit{),
then a necessary and sufficient condition that }$g=g^{\sigma }$\textit{\ for
some linear }$\sigma :Y\rightarrow \mathbb{R}$\textit{\ is the range
condition}%
\begin{equation*}
\mathcal{R}(K)\neq K(\mathcal{N}(\rho )).
\end{equation*}%
\textit{In this case }$\sigma $ \textit{is continuous provided }$K(\mathcal{N%
}(\gamma ))$\textit{\ is closed and complemented.}

\bigskip

\noindent \textbf{Proof. }Here $K|\mathcal{N}(\rho )$ is linear, and so $K(%
\mathcal{N}(\rho ))$ is a vector subspace of $Y.$

Suppose first that $\mathcal{R}(K)\neq K(\mathcal{N}(\rho )).$ Then there is 
$u\in X$ with $K(u)\notin K(\mathcal{N}(\rho )).$ So $K(u)\neq 0$ and $\rho
(u)\neq 0,$ so that $g(u)\neq 1.$ W.l.o.g. we may assume that $\rho (u)=1.$
Indeed, by Theorem 3.1, $K(u/\rho (u))=\lambda _{u}(\rho (u)^{-1})K(u),$ and
so $K(u/\rho (u))\notin K(\mathcal{N}(\rho ))$ by linearity of $K|\mathcal{N}%
(\rho ).$

Step 1. We first prove the result under the assumption that $Y=\langle
K(X)\rangle ,$ the span here being assumed a closed subspace.

We begin by defining a continuous linear map $\sigma $ by setting: 
\begin{equation}
\sigma (y):=\left\{ 
\begin{array}{cc}
0, & y\in K(\mathcal{N}(\rho )), \\ 
t(g(u)-1), & y=tK(u).%
\end{array}%
\right.  \tag{$\sigma _{A}$}
\end{equation}%
The two clauses are thus mutually exclsusive and so%
\begin{equation}
\sigma K(x)=g(x)-1  \tag{$Eq$}
\end{equation}%
certainly holds for the one vector $x=u.$

We first decompose $K$ into summands and likewise $g$ into factors, by
projecting along $\langle u\rangle .$ On these we act with $\sigma $, as $%
\sigma $ has non-zero effect only on the $K$-image $\langle u\rangle $.
Thereafter we reassemble the components.

As $\rho (x-\rho (x)u)=0$ and $g(x-\rho (x)u)=1,$%
\begin{eqnarray*}
K(x) &=&K([x-\rho (x)u]\circ _{\rho }\rho (x)u)=K(x-\rho (x)u)+1.K(\rho
(x)u): \\
K(x) &=&K(x-\rho (x)u)+K(\rho (x)u).
\end{eqnarray*}%
Since $x-\rho (x)u\in \mathcal{N}(\rho )$ and as $\sigma =0$ on $K(\mathcal{N%
}(\rho ))$ and is linear, applying $\sigma $ gives%
\begin{equation}
\sigma K(x)=\sigma K(\rho (x)u).  \tag{$A1$}
\end{equation}%
As $K(u)\neq 0$ and $g(u)\neq 1$ we may put (by Theorem 3.1)%
\begin{equation}
K(\rho (x)u)=\lambda _{w}(\rho (x))K(u).  \tag{$A2$}
\end{equation}%
Here $\lambda _{w}$ is defined by the formula of Theorem 3.1. By $(A1)$ and
applying $\sigma $ to $(A2),$%
\begin{equation}
\sigma K(x)=\sigma K(\rho (x)u)=\lambda _{w}(\rho (x))\sigma K(u)=\lambda
_{w}(\rho (x))[g(u)-1].  \tag{$A3$}
\end{equation}

This completes the action on the $K$ side.

We decompose $g$ similarly by $(M),$ as $\mathcal{N}(\rho )=\mathcal{N}%
(\gamma ):$ 
\begin{equation}
g(x)=g(x-\rho (x)u)\cdot g(\rho (x)u):\qquad g(x)=g(\rho (x)u).  \tag{$A4$}
\end{equation}%
Here again $\rho (x-\rho (x)u)=0,$ so $(x-\rho (x)u)\circ _{\rho }\rho
(x)u=x.$

We now act on the $g$ side.

By Lemma 6.2 on non-triviality, $g_{u}(tu)\neq 0$ for $t\neq 0,$ so we may
apply Lemma 6.1 (on switching). So%
\begin{eqnarray}
(g(\rho (x)u)-1)K(u) &=&K(\rho (x)u)[g(u)-1]  \notag \\
&=&[g(u)-1]\lambda _{w}(\rho (x))K(u)\qquad \text{(from }(A2)\text{),} 
\notag \\
(g(\rho (x)u)-1)[g(u)-1] &=&\lambda _{w}(\rho (x))[g(u)-1]^{2}\qquad \text{%
(applying }\sigma \text{):}  \notag
\end{eqnarray}%
\begin{equation}
(g(\rho (x)u)-1)=\lambda _{w}(\rho (x))[g(u)-1]\qquad \text{(cancelling),} 
\tag{$A5$}
\end{equation}%
as $g(u)-1\neq 0.$

We now reassemble the components. Combining, $(A5)$ with $(A4)$ and $(A3)$
gives%
\begin{equation*}
(g(x)-1)=\lambda _{w}(\rho (x))[g(u)-1]=\sigma (K(\rho (x)u))=\sigma K(x).
\end{equation*}

So $(Eq)$ holds for \textit{all} vectors $x\in X.$ This completes the
reassembly.

\bigskip

Step 2. If $Y\neq \langle K(X)\rangle ,$ choose in $Y$ a subspace $Z$
complementary to $\langle K(X)\rangle $ and define $\sigma $ as above on $%
\langle K(X)\rangle $; then extend by taking $\sigma =0$ on $Z.$

\bigskip

We turn to the converse and suppose now that $g=g^{\sigma }$ for some linear 
$\sigma :Y\rightarrow \mathbb{R}.$ We show that $(\sigma _{A})$ holds for
some $u\in X$, from which the range condition will follow. By $(GFE),$ $%
\sigma (y)=g(x)-1$ whenever $y=K(x).$ Further, as $(N_{A})$ holds, $g(x)=1$
for $x\in \mathcal{N}(\rho )$ and so, since%
\begin{equation*}
K(x)+K(x)=K(x)+g(x)K(x)=K(x\circ _{\rho }x)=K(x)+K(x)+\sigma (K(x))K(x),
\end{equation*}%
we conclude that%
\begin{equation*}
\sigma (K(x))=0,
\end{equation*}%
whether or not $K(x)=0.$ That is, $\sigma (y)=0$ for $y\in K(\mathcal{N}%
(\rho )),$ as in the first clause of $(\sigma _{A}).$

Since $g\neq 1$ there is $u\in X$ with $g(u)\neq 1.$ By $(N_{A}),$ $\rho
(u)\neq 0.$ The kernel $\mathcal{N}(\rho )$ is of co-dimension 1 in the
space $X,$ so that $X$ is the span of $u$ and $\mathcal{N}(\rho ).$ The
assumption $g=g^{\sigma }$ with $\sigma $ linear now gives for $t\neq 0$ 
\begin{equation*}
\sigma (K(u)=g(u)-1\neq 0:\qquad \sigma (tK(u))=t(g(u)-1),
\end{equation*}%
as in the second clause of $(\sigma _{A}).$ As the two clauses are
exclusive, $K(u)\notin K(\mathcal{N}(\rho )).$\hfill $\square $

\bigskip

Since, as above, the kernel $\mathcal{N}(\rho )$ is of co-dimension 1 in the
space $X,$ it seems natural to view the condition 
\begin{equation*}
\mathcal{R}(K)=K(\mathcal{N}(\rho ))
\end{equation*}%
as a type of degeneracy which we will refer to here (only) as yielding 
\textit{improper }solutions of\textit{\ }$(GFE).$

\bigskip

\textbf{An improper example} where $\mathcal{R}(K)=K(\mathcal{N}(\rho ))$ is
given by taking $X=Y=\mathbb{R}^{2},\rho (x)=x_{1}$ and $g(x)=1+x_{1}$ and $%
K(x)=(0,x_{2}).$ Then $(GFE)$ is satisfied, since%
\begin{equation*}
K(x\circ _{\rho }y)=(0,x_{2})+(1+x_{1})(0,y_{2}).
\end{equation*}%
Here $\mathcal{N}(\rho )=\{x:x_{1}=0\}=\mathcal{N}(\gamma ),$ so that $%
\mathcal{R}(K)=K(\mathcal{N}(\rho ))=\mathcal{N}(\rho ),$ and the
contradiction that $\sigma (0,x_{2})=x_{1}$ follows from the condition%
\begin{equation*}
x_{2}+(1+x_{1})y_{2}=x_{2}+(1+\sigma (0,x_{2}))y_{2}.
\end{equation*}

\bigskip

In the case $(N_{A})$ Prop. 7.1A above shows that all but the improper
solutions of $(GFE)$ are homomorphies. In the alternative case $(N_{B})$ all
solutions of $(GFE)$ are homomorphies, as we now show.

\bigskip

\noindent \textbf{Proposition 7.1B.}\textit{\ If }$(K,g)$\textit{\ satisfies 
}$(GFE$-$\rho )$\textit{\ and }$(N_{B}),$ \textit{that is,}%
\begin{equation*}
K(\mathcal{N}(\rho ))\subseteq \langle K(w)\rangle _{Y},
\end{equation*}%
\textit{for some} $w\in \mathcal{N}(\rho ),$ \textit{then }$g=g^{\sigma }$%
\textit{\ for some linear }$\sigma :Y\rightarrow \mathbb{R}$ \textit{which
is continuous, provided }$K(\mathcal{N}^{\ast }(\gamma ))$ \textit{is closed
complemented.}

\bigskip

\noindent \textbf{Proof. }Here $V_{0}:=\mathcal{N}^{\ast }(\gamma )$ $=%
\mathcal{N}(\gamma )\cap \mathcal{N}(\rho )$ is a subgroup of $\mathbb{G}%
_{\rho }(X),$ as%
\begin{equation*}
K(x+y)=K(x\circ _{\rho }y)=K(x)+g(x)K(y)=K(x)+K(y),
\end{equation*}%
and so $K|V_{0}$ is a homomorphism from $\mathbb{G}_{\rho }(X)$ to $Y$.
Since $V_{0}$ is a subspace of $\mathcal{N}(\rho ),$ we copy the argument of
Proposition 7.1A working with the linear map $K|V_{0}$ with $V_{0}\subseteq 
\mathcal{N}(\gamma )$ as a replacement for $K|\mathcal{N}(\rho );$ so if $%
g=g^{\sigma }$ is to hold, then $K|V_{0}:V_{0}\rightarrow \mathcal{N}(\sigma
).$

In $\mathcal{N}(\rho )$ choose a subspace $V_{1}$ complementary to $V_{0},$
and let $\pi _{i}:X\rightarrow V_{i}$ denote projection onto $V_{i}$. Notice
that for any $v\in \mathcal{N}(\rho ),$ as $\pi _{0}(v)\in \mathcal{N}(\rho
) $ and $\pi _{0}(v)\in \mathcal{N}(\gamma )$ 
\begin{eqnarray*}
K(v) &=&K(\pi _{0}(v)\circ \pi _{1}(v))=K(\pi _{0}(v))+g(\pi _{0}(v))K(\pi
_{1}(v)) \\
&=&K(\pi _{0}(v))+K(\pi _{1}(v)).
\end{eqnarray*}

Fix a non-zero $v_{1}\in V_{1}\subseteq \mathcal{N}(\rho );$ then $%
V_{1}=\langle v_{1}\rangle ,$ since by Lemma 2.1 $\gamma =\log g$ is linear
on $\mathcal{N}(\rho )$ and so $\mathcal{N}^{\ast }(\gamma )$ either equals $%
\mathcal{N}(\rho )$ or has co-dimension 1 in $\mathcal{N}(\rho )$. We can
see this directly as follows. Since $v_{1}\notin V_{0},$ $\gamma (v_{1})\neq
0,$ so replacing $v_{1}$ by $v_{1}/\gamma (v_{1}),$ w.l.o.g. $\gamma
(v_{1})=1$. For any $z\in \mathcal{N}(\rho ),$ $z-\gamma (z)v_{1}\in 
\mathcal{N}(\gamma ),$ as 
\begin{equation*}
\gamma (z-\gamma (z)v_{1})=\gamma (z)-\gamma (z)=0.
\end{equation*}%
Likewise, for such $z,$ as $\rho (v_{1})=0,$%
\begin{equation*}
\rho (z-\gamma (z)v_{1})=\rho (z)-\gamma (z)\rho (v_{1})=0-0=0,
\end{equation*}%
i.e. $v_{0}:=z-\gamma (z)v_{1}\in V_{0}$ and so%
\begin{equation*}
z=v_{0}+\gamma (z)v_{1}T\text{ i.e. }\mathcal{N}(\rho )=V_{0}+\langle
v_{1}\rangle .
\end{equation*}

If $\rho $ is not identically zero, again fix $u\in X$ with $\rho (u)=1.$
Then $x\mapsto \pi _{u}(x)=x-\rho (x)u$ is again (linear) projection onto $%
\mathcal{N}(\rho )$. If $\rho \equiv 0,$ set $u$ below to $0.$ Whether or
not $\rho \equiv 0,$ as $\rho (x-\rho (x)u)=0$ take $z:=x-\rho (x)u\in 
\mathcal{N}(\rho );$ then for some $v_{0}\in V_{0}$ and some scalar $\alpha $%
\begin{equation*}
x=v_{0}+\alpha v_{1}+\rho (x)u=v\circ _{\rho }\rho (x)u.
\end{equation*}%
So w.l.o.g. provided $K(v_{1})\neq 0\neq K(u)$%
\begin{equation*}
Y=\langle K(V_{0}),K(v_{1}),K(u)\rangle .
\end{equation*}

We first show that these \textquotedblleft generators\textquotedblright\ are
distinct. Recall that $g(v_{1})\neq 1$ as $v_{1}\notin V_{0}$ and that for
some $w$ with $K(w)\neq 0$ 
\begin{equation*}
K(\mathcal{N}(\rho ))=K(V_{0}+\langle v_{1}\rangle )\subseteq \langle
K(w)\rangle _{Y}.
\end{equation*}

Suppose first that $K(u)=K(v)$ for some $v\in V_{0}+\langle v_{1}\rangle .$
Then, as $K(u)=K(v)\in \langle K(w)\rangle _{Y}$ and $K(w)\neq 0,$%
\begin{equation*}
K(u)=\lambda _{w}(u)K(w)=\lambda _{w}(v)K(w).
\end{equation*}%
So, since $\lambda _{w}$ is montonic,%
\begin{equation*}
u=v\in \mathcal{N}(\rho ),
\end{equation*}%
contradicting that $\rho (u)=1.$

Next suppose that $K(v_{1})=K(v_{0}),$ for some $v_{0}\in V_{0}.$ Then,
since $-v_{0}+v_{1}=-v_{0}\circ _{\rho }v_{1}$ and $g(v_{0})=1,$ 
\begin{eqnarray*}
0 &=&-K(v_{0})+g(v_{0})K(v_{1})=K(-v_{0}+v_{1}) \\
&=&K(v_{1})-g(v_{1})K(v_{0})=K(v_{1})-g(v_{1})K(v_{1}) \\
&=&(1-g(v_{1}))K(v_{1})=0.
\end{eqnarray*}%
So $K(v_{1})=0,$ a contradiction.

So the following defines a continuous linear map $\sigma :Y\rightarrow 
\mathbb{R}$ : 
\begin{equation}
\sigma (y)=\left\{ 
\begin{array}{cc}
0, & y\in K(V_{0}), \\ 
t(g(v_{1})-1) & y=tK(v_{1}), \\ 
t(g(u)-1), & y=tK(u).%
\end{array}%
\right.  \tag{$\sigma _{B}$}
\end{equation}%
So $(Eq)$ holds for the two vectors $x=v_{1}$ and $x=u.$

As with $(A1)$ in Proposition 7.1A, via Lemma 6.1 (on switching),%
\begin{equation}
\sigma K(\rho (x)u)=g(\rho (x)u)-1,  \tag{$B1$}
\end{equation}%
\begin{equation}
\sigma K(\alpha v_{1})=g(\alpha v_{1})-1.  \tag{$B2$}
\end{equation}%
Since $v_{i}$ are in $\mathcal{N}(\rho ),$ 
\begin{eqnarray*}
K(x) &=&K(v_{0}+\alpha v_{1}+\rho (x)u) \\
&=&K(v_{0})+K(\alpha v_{1})+g(\alpha v_{1})K(\rho (x)u).
\end{eqnarray*}%
Using $(\sigma _{B})$ and applying $\sigma $ gives%
\begin{eqnarray*}
\sigma K(x) &=&0+[g(\alpha v_{1})-1]+g(\alpha v_{1})[g(\rho (x)u)-1]\text{
\qquad (by }(B2)\text{ and }(B1)\text{)} \\
&=&[g(\alpha v_{1})-1]-g(\alpha v_{1})+g(\alpha v_{1})g(\rho (x)u) \\
&=&g(v_{0})g(\alpha v_{1})g(\rho (x)u)-1 \\
&=&g(x)-1.
\end{eqnarray*}

If $Y\neq \langle K(V_{0}),K(v_{1}),K(u)\rangle ,$ this span being assumed a
closed subspace, choose in $Y$ a subspace $Z$ complementary to $\langle
K(V_{0}),K(v_{1}),K(u)\rangle ,$ and define $\sigma $ as above on $\langle
K(V_{0}),K(v_{1}),K(u)\rangle ;$ then extend by taking $\sigma =0$ on $Z.$%
\hfill $\square $

\bigskip

In Lemma 7.1 below we refer to the defining equation $(g^{\sigma })$ in \S 1.

\bigskip

\noindent \textbf{Lemma 7.1. }\textit{If }$(K,g)$\textit{\ satisfies }$(GFE$-%
$\rho )$\textit{\ non-trivially, then }$g$\textit{\ is uniquely determined
by }$K.$ \textit{In particular, if }$g^{\sigma }=g=$\textit{\ }$g^{\tau }$%
\textit{, then }$\sigma =\tau $\textit{\ on} $K(X).$ \textit{Furthermore,}%
\begin{equation*}
\sigma (K(u\circ _{\rho }v))=\sigma (K(u))\circ _{\iota }\sigma (K(v))\qquad
(\iota (t)\equiv t).
\end{equation*}

\noindent \textbf{Proof. }For given $K,$ suppose both $(K,g)$ and $(K,h)$
satisfy $(GFE)$. As $K$ is non-trivial, we fix $v\in X$ with $K(v)\neq 0$.
Then for arbitrary $u\in X$%
\begin{equation*}
K(u)+h(u)K(v)=K(u\circ _{\rho }v)=K(u)+g(u)K(v),
\end{equation*}%
so $g(u)=h(u).$ So if $g^{\tau }=g=g^{\sigma },$ then 
\begin{equation*}
\sigma (K(u))=g^{\sigma }(u)-1=g^{\tau }(u)-1=\tau (K(u)),
\end{equation*}%
as claimed. Checking the final assertion is routine: see \S 9.5 (Appendix).
\hfill $\square $

\bigskip

We may now pass to the key existence theorem.

\bigskip

\noindent \textbf{Theorem 7.2. }\textit{If }$(K,g)$\textit{\ satisfies }$%
(GFE $-$\rho ),$\textit{\ then, unless} $\mathcal{R}(K)=K(\mathcal{N}(\rho
)),$\textit{\ there is a unique linear map} $\sigma :Y\rightarrow \mathbb{R}$
\textit{such that}%
\begin{equation*}
\sigma (K(x))+1=g(x)\qquad (x\in X).
\end{equation*}%
\textit{The map }$\sigma $\textit{\ is continuous, provided }$K$\textit{\
has closed complemented range.}

\bigskip

\noindent \textbf{Proof. }By Corollary 6.2, one of $(N_{A}^{=})$ or $%
(N_{B}^{=})$ holds, and so either Proposition 7.1A or 7.1B implies the
existence of $\sigma ,$ and its continuity conditional on $K$ having closed
range. Its uniqueness is assured by Lemma 7.1. \hfill $\square $

\bigskip

\noindent \textbf{8. The generalized Goldie equation.}

This section is devoted to demonstrating in Theorem 8.1 below that $(GGE)$
is reducible to $(GFE)$. Our main tool is Theorem 3.1, and we will also use
Theorem 2.1 (the Index Theorem). For $(GGE)\ $to conform with our study of $%
(GFE$-$\rho ),$ we assume here that, just like $\eta _{\rho },$ the
non-negative inner auxiliary $h:X\rightarrow \lbrack 0,\infty )$ preserves
positivity on $\mathbb{G}_{h}^{+}(X):=\{x\in X:h(x)>0\}:$%
\begin{equation*}
h(u+h(u)v)>0\text{ for }u,v\in \mathbb{G}_{h}^{+}(X).
\end{equation*}%
This assumption prompts the question of for which continuous functions $%
h:X\rightarrow \lbrack 0,\infty )$ does the binary operation $x\circ
_{h}y:=x+h(x)y$ \textit{preserve positivity}, i.e.%
\begin{equation*}
h(x),h(y)>0\Longrightarrow h(x+h(x)y)>0.
\end{equation*}%
It emerges that strengthening $\Longrightarrow $ above to $%
\Longleftrightarrow $ yields Chudziak's theorem, that $h$ satisfies $(GS)$
(i.e. $h=\eta _{\rho }$ for some $\rho $). For details see [Chu2] (cf.
[Chu4]). Our hypothesis is thus weaker; however, this preservation combined
with $(GGE)$ yields some similar connections with $(GS)$ below.

We will need to know the connection between the null spaces of the inner and
outer auxiliaries. Recall that $K(0)=0.$

\bigskip

\noindent \textbf{Lemma 8.1. }\textit{If }$(K,g,h)$ \textit{satisfies }$%
(GGE) $\textit{\ and }$K(w)\neq 0$ \textit{for some }$w$\textit{, then}%
\begin{equation*}
g(x)=0\Longleftrightarrow h(x)=0\qquad (x\in X),
\end{equation*}%
\textit{\ so that }%
\begin{equation*}
\mathbb{G}_{h}^{+}(X):=\{x\in X:h(x)>0\}=\mathbb{G}^{+}(X)=\mathbb{G}%
_{g}^{+}(X):=\{x\in X:g(x)>0\}.
\end{equation*}

\noindent \textbf{Proof. }If $h(a)=0,$ then $g(a)=0,$ since $K(w)\neq 0$ and%
\begin{equation*}
K(a)=K(a+h(a)w)=K(a)+g(a)K(w).
\end{equation*}%
If one had $g(a)=0$ but $h(a)\neq 0,$ then, for any $x,$ taking $%
b:=h(a)^{-1}(x-a)$ gives%
\begin{equation*}
K(a)=K(a)+g(a)K(b)=K(a+h(a)b)=K(x).
\end{equation*}%
Hence $K$ is constant. But $K(0)=0,$ so $K(w)=0,$ a contradiction. \hfill $%
\square $

\bigskip

Our first result identifies a known partially `pexiderized' variant of the Go%
\l \k{a}b-Schinzel equation, $(PGS)$ below, studied in [Chu1], [Jab1]; see
[Jab3] for a fully pexiderized equation (cf. [Jab2]). The gist of the matter
is in Prop. 4.4 and the solution is given by%
\begin{equation*}
h_{u}(t):=\left\{ 
\begin{array}{ccc}
1+rt, & \text{for }r\neq 0,\text{ and then: } & g_{u}=g_{r,\theta
}=(1+rt)^{\theta /r}; \\ 
1, & \text{for }r=0,\text{ and then: } & g_{u}=g_{0,\theta }=e^{\theta t}.%
\end{array}%
\right.
\end{equation*}%
Here $g_{u}(t)=g(tu)$ and $h_{u}(t)=h(tu).$ The next Proposition links the
behaviour of the two auxiliary functions $g,h$ via the linking function $%
\lambda _{w}.$ Motivated by the notation $\langle u\rangle _{\rho }$ in Popa
groups, but now in the context of a \ Javor group, we write for $u\in X$%
\begin{equation*}
\langle u\rangle _{h}=\mathbb{G}_{h}^{+}(X)\cap \mathrm{Lin}\{u\}.
\end{equation*}

\bigskip

\noindent \textbf{Proposition 8.1. }\textit{For }$(K,h,g)$\textit{\
satisfying }$(GGE)$\textit{\ and for each }$w\in X$\textit{\ with }$K(w)\neq
0,$ \textit{and} $u\in \mathbb{G}^{+}(X),$%
\begin{equation*}
\lambda _{w}\left( \frac{h(a+h(a)b)}{h(a)h(b)}\right) =\frac{g(a+h(a)b)}{%
g(a)g(b)}\qquad (a,b\in \mathbb{G}_{h}^{+}(X)).
\end{equation*}%
\textit{In particular, if the auxiliary }$h_{u}$\textit{\ satisfies the Go\l 
\k{a}b-Schinzel equation, then }$g_{u}$\textit{\ satisfies a partially
pexiderized Go\l \k{a}b-Schinzel equation:}%
\begin{equation}
g(a+h(a)b)=g(a)g(b)\text{ }\qquad (a,b\in \langle u\rangle _{h}). 
\tag{$PGS$}
\end{equation}%
\textit{So }$g=g_{\gamma (u),\rho (u)}$\textit{\ (for some appropriate
parameters), and conversely if }$g$\textit{\ has this form, then }$g$\textit{%
\ satisfies }$(PGS)$ \textit{for }$h_{u}=1+\rho (u).$

\bigskip

\noindent \textbf{Proof. }We approach the action of $K$ on%
\begin{equation*}
a+h(a)b+h(a+h(a)b)h(a)h(b)w
\end{equation*}%
in two ways. For the approach to be valid we need $h(a+h(a)b)>0,$ which
comes from the assumed preservation of positivity (cf. Lemma 8.1). We
consider the two sides of the equality%
\begin{equation*}
K(a+h(a)b+h(a+h(a)b)h(a)h(b)w)=K(a+h(a)b)+g(a+h(a)b)K(h(a)h(b)w).
\end{equation*}%
Here, with $LHS$ for left-hand side etc.,%
\begin{eqnarray*}
LHS &=&K(a+h(a)[b+h(a+h(a)b)h(b)w]) \\
&=&K(a)+g(a)K(b+h(b)h(a+h(a)b)w) \\
&=&K(a)+g(a)[K(b)+g(b)K(h(a+h(a)b)w) \\
&=&K(a)+g(a)K(b)+g(a)g(b)K(h(a+h(a)b)w); \\
RHS &=&K(a)+g(a)K(b)+g(a+h(a)b)K(h(a)h(b)w).
\end{eqnarray*}%
Cancelling common terms on the two sides gives, in view of $g(a)g(b)\neq 0,$
that%
\begin{eqnarray*}
g(a)g(b)K(h(a+h(a)b)w) &=&g(a+h(a)b)K(h(a)h(b)w): \\
K\left( \frac{h(a+h(a)b)}{h(a)h(b)}w\right) &=&\frac{g(a+h(a)b)}{g(a)g(b)}%
K(w),\text{ }
\end{eqnarray*}%
on replacing $w$ appropriately (since $h(a)h(b)\neq 0$). Now for $K(w)\neq
0, $ apply Theorem 3.1. So if $h$ satisfies $(GS),$ then $g$ satisfies $%
(PGS).$ \hfill $\square $

\bigskip

\noindent \textbf{Corollary 8.1}. \textit{If }$(K,h,g)$\textit{\ satisfies }$%
(GGE),$\textit{\ then either }$K$\textit{\ is linear and }$g|\langle
w\rangle _{h}=h|\langle w\rangle _{h}$ \textit{for each }$w,$\textit{\ or }$%
h $\textit{\ satisfies }$(GS)$\textit{\ and }$g$\textit{\ satisfies }$(PGS).$

\bigskip

Before deducing this result we need to characterize the `contour behaviour'
of the link functions $\lambda _{u}(t)$ in response to parameter changes.

We recall that when the kernel function is non-zero at $u$ (i.e. $K(u)\neq 0$%
) the link function $\lambda _{u}(t)$ is either the identity function 
\textrm{id}$(t)=t,$ or for some $r>0$ $\lambda _{u}(t)=\varphi (rt)/\varphi
(r),$ where $\varphi (x)$ takes one of three \textit{contour types}, the
exponential, logarithmic, or power-$c$, all with domain parameter $r:$%
\begin{equation*}
\begin{array}{cc}
\varphi (x):=e^{x}-1: & \lambda _{u}(t)=(e^{rt}-1)/(e^{r}-1), \\ 
\varphi (x):=1+\log x\text{ }: & \lambda _{u}(t)=\log (1+rt)/\log (1+r), \\ 
\varphi (x):=(1+x)^{c}-1: & \lambda _{u}(t)=[(1+rt)^{c}-1]/[(1+r)^{c}-1].%
\end{array}%
\end{equation*}%
These are strictly monotone in $t$ and either convex or concave shaped. In
the power-$c$ type convexity arises for $c>1$ and (like the exponential) is
separated from its concave inverse function by the linear variant $\lambda
_{u}(t)=t$ arising from $c=1.$

Indeed, it is enough to note the relevant second derivative, $\varphi
^{\prime \prime }(rt),$ which according to type is%
\begin{equation*}
r^{2}e^{rt},\qquad -r^{2}(1+rt)^{-2},\qquad r^{2}c(c-1)(1+rt)^{c-2}.
\end{equation*}%
\bigskip

\noindent \textbf{Lemma 8.2.} \textit{For any point }$(x,y)$\textit{\ in the
positive quadrant other than }$(1,1),$\textit{\ there is at most one curve }$%
\lambda _{u}$\textit{\ in any of the three contour types with }$\lambda
_{u}(x)=y.$

\bigskip

\noindent \textbf{Proof. }For $s>0$ and $u$ with $K(u)\neq 0$ as above,
recall Cor. 4.3(i):%
\begin{equation*}
\lambda _{su}(t)=\lambda _{u}(st)/\lambda _{u}(s).
\end{equation*}%
From here it follows that scaling the vector $u$ by $s>0$ does not alter the 
\textit{contour type }of the link function $\lambda _{su}$ but merely scales
its domain parameter from $r$ to $sr.$ Indeed, given the tabulation above,
this follows from: 
\begin{equation*}
\lambda _{u}(st)/\lambda _{u}(s)=\left. \frac{\varphi (rst)}{\varphi (r)}%
\right/ \frac{\varphi (rs)}{\varphi (r)}=\frac{\varphi (rst)}{\varphi (rs)}.
\end{equation*}

The contours all have $(1,1)$ as fixed point, but on each side of $t=1$ the
curves of any one type are strictly \textit{monotone} in the domain
parameter $r$, as some routine calculus readily shows (Lemma 8.4 below). For
example, the exponential type curves decrease with $r$ to the left of $t=1$
and increase with $r$ to the right of $t=1.$

Hence for any point $(x,y)$ in the positive quadrant other than $(1,1)$
there is at most one curve $\lambda _{u}$ in each type with $\lambda
_{u}(x)=y.$ \hfill $\square $

\bigskip

\noindent \textbf{Proof of Corollary 8.1. }We apply Lemma 8.2 and consider $%
w $ with $K(w)\neq 0$. By Proposition 8.1, for all $s>0$ and $a,b\in \mathbb{%
G}_{h}^{+}(X)=\{x:h(x)>0\},$ since $K(sw)=\lambda _{w}(s)K(w)\neq 0,$%
\begin{equation*}
\lambda _{sw}\left( \frac{h(a+h(a)b)}{h(a)h(b)}\right) =\frac{g(a+h(a)b)}{%
g(a)g(b)}=\lambda _{w}\left( \frac{h(a+h(a)b)}{h(a)h(b)}\right) :
\end{equation*}%
\begin{equation}
\lambda _{sw}\left( \frac{h(a+h(a)b)}{h(a)h(b)}\right) =\lambda _{w}\left( 
\frac{h(a+h(a)b)}{h(a)h(b)}\right) .  \tag{$\ast $}
\end{equation}%
In the case when $\lambda _{w}(t)\equiv t$ this last equation holds, no
matter the value of $s$, since also $\lambda _{sw}(t)\equiv t.$ But that is
a very special case, which leads to a linear $K.$ To begin with, $K$ is
homogeneous, i.e. for all $t>0$ and all $w$ in $\mathbb{G}_{+}^{h}(X)$%
\begin{equation*}
K(tw)=tK(w).
\end{equation*}%
But in this case $g=h$ on $\langle w\rangle _{h}.$ Indeed, for $s,t>0,$ 
\begin{equation*}
sK(w)+g_{w}(s)tK(w)=K(sw+h_{w}(s)tw)=(s+h_{w}(s)t)K(w).
\end{equation*}%
This implies that $K$ is additive and so linear (by homogeneity):%
\begin{equation*}
K(a+b)=K\left( a+h(a)\frac{b}{h(a)}\right) =K(a)+\frac{g(a)}{h(a)}%
K(b)=K(a)+K(b).
\end{equation*}%
Here $h$ may be arbitrary with $g=h$ on each ray $\langle w\rangle _{h}.$

So suppose now that for some $w$ we have $\lambda _{w}\neq \mathrm{id.}$
Then for $s>0$ all the curves $\lambda _{sw}$ are of one contour type
differing only in their domain parameter. So $(\ast )$ contradicts Lemma 8.2
in that there may be at most one contour in any contour type passing through
a point unless that point is $(1,1).$ Hence%
\begin{equation*}
\frac{h(a+h(a)b)}{h(a)h(b)}=1=\frac{g(a+h(a)b)}{g(a)g(b)}.
\end{equation*}%
It now follows that $h$ satisfies the $(GS)$ equation and $g$ satisfies the
pexiderized variant $(PGS).$\hfill $\square $

\bigskip

We need a further (folk-lore) result\footnote{%
Recalled by Prof Chudziak at the 20th ICFE.}.

\bigskip

\noindent \textbf{Lemma 8.3.} \textit{For }$\rho $\textit{\ homogeneous on }$%
\mathbb{G}_{\rho }(X),$\textit{\ if }$\eta _{\rho }(u)=1+\rho (u)$\textit{\
satisfies }$(GS)$\textit{, then }$\rho $\textit{\ is linear.}

\bigskip

\noindent \textbf{Proof.} Fix $u,v$ $\in \mathbb{G}_{h}^{+}(X)$ and consider 
$\alpha ,\beta $ with $t=:\eta _{\rho }(\alpha u)=1+\rho (\alpha u)>0$ and $%
\eta _{\rho }(\beta u)>0.$ Then 
\begin{eqnarray*}
\rho (\alpha u+\beta v) &=&\eta _{\rho }(\alpha u+\eta _{\rho }(\alpha
u)\beta v/t)-1=(1+\rho (\alpha u))(1+\rho (\beta v/t))-1 \\
&=&\rho (\alpha u)+(1+\rho (\alpha u))\rho (\beta v/t)=\alpha \rho (u)+t\rho
(\beta v/t) \\
&=&\alpha \rho (u)+\beta \rho (v).\qquad \qquad \qquad \square
\end{eqnarray*}%
\textbf{\ }

\bigskip

\noindent \textbf{Theorem 8.1. }\textit{Suppose }$X$ \textit{is a normed
vector space and }$(K,h,g)$ \textit{satisfies }$(GGE)$ \textit{with }$K$ 
\textit{non-trivial }(\textit{i.e. there is }$w\in X$ \textit{with }$%
K(w)\neq 0\mathrm{).}$ \textit{Then either }$K$ \textit{is linear or else,
for some continuous linear }$\rho ,$%
\begin{equation*}
h(su)=1+s\rho (u)\qquad (s>0,u\notin \mathcal{N}(K)).
\end{equation*}%
\textit{In particular, }$(K,g)$ \textit{satisfies }$(GFE$-$\rho )$\textit{\
and so, for any }$u$ \textit{with }$\rho (u)=1,$%
\begin{equation*}
g(x)=e^{\alpha (x)}(1+\rho (x))^{\beta }\qquad (x\in \mathbb{G}_{\rho }(X)),
\end{equation*}%
\textit{where} 
\begin{equation*}
\alpha (x):=\gamma (x-\rho (x)u)\quad (x\in \mathbb{G}_{\rho }(X))\text{ }
\end{equation*}%
\textit{\ is linear and }$\alpha (u)=0$ \textit{and }$\beta =\gamma (u)/\log
2.$ \textit{Thus there are four cases:}%
\begin{equation*}
\begin{array}{ccc}
h_{u}(s)=g_{u}(s)=1, &  & \rho (u)=\gamma (u)=0, \\ 
h(su)=1,\qquad g(su)=e^{s\gamma (u)}, &  & \rho (u)=0\text{ and\textit{\ }}%
\gamma (u)\neq 0, \\ 
h_{u}(s)=1+s\rho (u),\qquad g_{u}(s)=1, &  & \rho (u)\neq 0\text{ and}%
\mathit{\ }\gamma (u)=0, \\ 
h_{u}(s)=1+s\rho (u),\qquad g_{u}(s)=(1+s\rho (u))^{(\gamma /\rho )}, &  & 
\gamma (u)\neq 0\neq \rho (u).%
\end{array}%
\end{equation*}

\textit{The form of }$K$\textit{\ may be read off from }$[$\textit{BinO8,Th
4A, 4B}$]$\textit{.}

\bigskip

\noindent \textbf{Proof. }By Corollary 5.1, $h_{u}(s)=1+s\rho
(u)=h(su)=h_{su}(1)=1+\rho (su),$ so that $\rho (u)$ is homogeneous. By
Lemma 5.3, $\rho $ is linear and by assumption continuous. Hence the
equation $(GGE)$ has the form $(GFE),$ i.e.%
\begin{equation*}
K(u\circ _{\rho }v)=K(u)+g(u)K(v).
\end{equation*}%
The remaining assertions follow from Theorem 2.1 (the Index Theorem).\hfill $%
\square $

\bigskip

The four cases above can also be reached from Theorem 3.1 by four direct but
laborious computations using Prop. 4.4: see \S 9.6 (Appendix). Theorem 7.2
puts these last conclusions into perspective, since $(GFE)$ above reduces to
a homomorphism between $\mathbb{G}_{\rho }(X)$ and $\mathbb{G}_{\sigma }(Y)$
for some $\sigma ,$ unless the range condition is violated. Recall from
[BinO8, Th.2] that, being abelian, $K(\mathcal{N}(\rho ))$ is either
included in $\mathcal{N}(\sigma )$ (the $N_{A}$-case of Theorem 7.1 with $%
\mathcal{N}(\rho )\subseteq \mathcal{N}(\gamma )$) or lies along a radius in 
$Y$ (the $N_{B}$-case). In the first case, by [BinO8,Th 4A], $K$ exhibits
linear and power types of behaviour (the latter only if $\rho (u)=1$ for
some $u,$ the behaviour becoming logarithmic in the limit when $\sigma
(K(u))=0$). Otherwise, by [BinO8, Th 4B], $K$ exhibits exponential and power
type behaviour (the latter again if $\rho (u)=1$ for some $u,$ again
becoming logarithmic when $\sigma (K(u))=0$).

\bigskip

We close by verifying how the $\lambda _{u}(t;r)$ contour rises or falls as
the domain parameter $r$ rises, according to contour type, and according to
which side of of $t=1;$ it is here at $t=1$ where behaviour reverses (cf. $%
\lambda _{u}(t)\simeq e^{r(t-1)},$ using a large $r$ approximation). The
calculations check for monotonicity with a simple scheme based on the form $%
\lambda (t)=\varphi (rt)/\varphi (r).$

\bigskip

\noindent \textbf{Lemma 8.4.} (a) \textit{Consider }$\lambda
_{u}(t;r)=\varphi (rt)/\varphi (r)$ \textit{with }$\varphi (x)=e^{x}-1.$%
\textit{\ Then:}

\noindent i) \textit{for }$0<t<1$\textit{\ the }$\lambda $ \textit{contours
fall as }$r$\textit{\ rises: if }$0<r<R,$\textit{\ then}

\begin{equation*}
\lambda (t;R)<\lambda (t;r);
\end{equation*}%
\noindent ii) \textit{by reciprocation, at any }$s>1,$\textit{\ if }$0<r<R,$%
\textit{\ then}%
\begin{equation*}
\lambda (s;r)<\lambda (t;R);
\end{equation*}%
\noindent iii) \textit{by inversion, the results in (i) and (ii) are
reversed for }$\varphi (x)=\log (1+x).$

(b) \textit{Consider }$\lambda _{u}(t;r)=\varphi (rt)/\varphi (r)$ \textit{%
with }$\varphi (x)=(1+x)^{c}-1$ \textit{and }$c<1.$ \textit{Then:}

\noindent i) \textit{for }$0<t<1$\textit{\ the }$\lambda $ \textit{contours
fall as }$r$ \textit{rises: if }$0<r<R,$\textit{\ then}

\begin{equation*}
\lambda (t;R)<\lambda (t;r);
\end{equation*}%
\noindent ii) \textit{by reciprocation, at any }$s>1,$\textit{\ if }$0<r<R,$%
\textit{\ then}%
\begin{equation*}
\lambda (s;r)<\lambda (t;R);
\end{equation*}%
\noindent iii) \textit{by inversion, the results in (i) and (ii) are
reversed for }$c>1.$

\bigskip

\noindent \textbf{Proof.} We begin by computing that%
\begin{eqnarray*}
\frac{\partial }{\partial r}\left( \frac{\varphi (rt)}{\varphi (r)}\right)
&=&\frac{t\varphi ^{\prime }(rt)\varphi (r)-\varphi (rt)\varphi ^{\prime }(r)%
}{\varphi (r)^{2}}=\frac{\varphi (rt)\varphi (r)}{\varphi (r)^{2}}\left( t%
\frac{\varphi ^{\prime }(rt)}{\varphi (rt)}-\frac{\varphi ^{\prime }(r)}{%
\varphi (r)}\right) , \\
\frac{\partial }{\partial r}\left( \frac{\varphi ^{\prime }(rt)}{\varphi (rt)%
}\right) &=&\frac{t\varphi ^{\prime \prime }(rt)\varphi (rt)-t\varphi
^{\prime }(rt)\varphi ^{\prime }(rt)}{\varphi (rt)^{2}}.
\end{eqnarray*}%
(a) (i) As $\varphi ^{\prime }=\varphi ^{\prime \prime }=e^{x},$ we have for 
$0<t<1$%
\begin{equation*}
te^{rt}(e^{rt}-1)-te^{rt}e^{rt}=-te^{rt}<0\text{:\qquad }\frac{\partial }{%
\partial r}\left( t\frac{\varphi ^{\prime }(rt)}{\varphi (rt)}\right) <0.
\end{equation*}%
Since $rt<r$ for $0<t<1,$ we have%
\begin{equation*}
t\frac{\varphi ^{\prime }(rt)}{\varphi (rt)}<\frac{\varphi ^{\prime }(r)}{%
\varphi (r)},\text{ so }\frac{\partial }{\partial r}\left( \frac{\varphi (rt)%
}{\varphi (r)}\right) <0:\qquad \frac{\varphi (rt)}{\varphi (r)}\text{ is
decreasing in }r.
\end{equation*}%
Hence, if $0<r<R$ and $0<t<1,$ then%
\begin{equation*}
\varphi (rt)/\varphi (r)>\varphi (Rt)/\varphi (R).
\end{equation*}

(ii) With $0<t<1,$ note that $\varphi (r)/\varphi (rt)$ is increasing in $%
r>0.$ Writing $\rho =rt$ and $s=1/t>1,$ we see that $\varphi (\rho
s)/\varphi (\rho )$ is increasing in $\rho $ ($t$ being fixed.)

(iii) The final assertion follows because $y=\varphi (x)=e^{x}-1$ has
inverse $x=\log (1+y).$

(b) Here $\varphi ^{\prime }=c(1+x)^{c-1}$ and $\varphi ^{\prime \prime
}=c(c-1)(1+x)^{c-2},$ so for $c,t<1$%
\begin{equation*}
t\varphi ^{\prime \prime }(rt)\varphi (rt)-\varphi ^{\prime }(rt)\varphi
^{\prime }(rt)=tc(c-1)(1+rt)^{c-2}[(1+rt)^{c}-1]-c^{2}(1+rt)^{2c-2}<0.
\end{equation*}%
since $(1+rt)>1.$ This leads to (i) and (ii), exactly as in (a).

(iii) Since $y=\varphi _{c}(x)=(1+x)^{c}-1$ has as inverse $x=\varphi
_{1/c}(y)=(1+y)^{1/c}-1,$ the assertion (iii) now follows from (i) and (ii)
with reversal. \hfill $\square $

\bigskip

\textbf{References}

\noindent \lbrack Acz] J. Acz\'{e}l. \textsl{Lectures on Functional
Equations and Their Applications.} Academic Press, New York,1966.

\noindent \lbrack AczD] J. Acz\'{e}l and J. Dhombres, \textsl{Functional
equations in several variables.} Encycl. Math. App.\textbf{\ 31}, Cambridge
University Press, 1989.\newline
\noindent \lbrack AczG] J. Acz\'{e}l and S. Go\l \k{a}b, Remarks on
one-parameter subsemigroups of the affine group and their homo- and
isomorphisms. \textsl{Aequat. Math.} \textbf{4} (1970), 1--10.\newline
\noindent \lbrack BinGT] N. H. Bingham, C. M. Goldie and J. L. Teugels, 
\textsl{Regular variation}, 2nd ed., Cambridge University Press, 1989 (1st
ed. 1987).\newline
\noindent \lbrack BinO1] N. H. Bingham and A. J. Ostaszewski, Homotopy and
the Kestelman-Borwein-Ditor theorem. \textsl{Canad. Math. Bull.} \textbf{54}
(2011), 12--20. \newline
\noindent \lbrack BinO2] N. H. Bingham and A. J. Ostaszewski, Beurling slow
and regular variation. \textsl{Trans. London Math. Soc. }\textbf{1} (2014),
29-56.\newline
\noindent \lbrack BinO3] N. H. Bingham and A. J. Ostaszewski, Cauchy's
functional equation and extensions: Goldie's equation and inequality, the Go%
\l \k{a}b-Schinzel equation and Beurling's equation. \textsl{Aequat. Math.} 
\textbf{89} (2015), 1293--1310.\newline
\noindent \lbrack BinO4] N. H. Bingham and A. J. Ostaszewski, Beurling
moving averages and approximate homomorphisms. \textsl{Indag. Math. }\textbf{%
27} (2016), 601-633 (fuller version: arXiv1407.4093).\newline
\noindent \lbrack BinO5] {N. H. Bingham and A. J. Ostaszewski, }General
regular variation, Popa groups and quantifier weakening. \textsl{J. Math.
Anal. Appl.} \textbf{483} (2020) 123610, 31 pp. (arXiv1901.05996). \newline
\noindent \lbrack BinO6] {N. H. Bingham and A. J. Ostaszewski, }Sequential
regular variation: extensions of Kendall's theorem, \textsl{Quart. J. Math. }%
\textbf{71} (2020), no. 4, 1171--1200 (arXiv:1901.07060).\newline
\noindent \lbrack BinO7] {N. H. Bingham and A. J. Ostaszewski, Extremes and
regular variation. \textsl{A lifetime of excursions through random walks and
L\'{e}vy processes}, 121--137, \textsl{Progr. Probab.} \textbf{78}, Birkh%
\"{a}user/Springer, Cham, 2021 (arXiv2001.05420).}\newline
\noindent \lbrack BinO8] {N. H. Bingham and A. J. Ostaszewski, }%
Homomorphisms from Functional Equations: II. The Goldie Equation\textbf{. (}%
arXiv:1910.05816; originally under the title: Multivariate general regular
variation: Popa groups on vector spaces)\newline
\noindent \lbrack BinO9] {N. H. Bingham and A. J. Ostaszewski, }The Go\l 
\k{a}b-Schinzel and Goldie functional equations in Banach algebras,
arXiv:2105.07794.\newline
\noindent \lbrack BinO10] {N. H. Bingham and A. J. Ostaszewski, }\textsl{%
Category and measure: Infinite combinatorics, topology and Groups.}
Cambridge University Press, Cambridge, 2024, forthcoming.\newline
\noindent \lbrack BinO11] {N. H. Bingham and A. J. Ostaszewski.}
Parthasarathy, shift-compactness and infinite combinatorics, \textsl{Indian
J Pure and Applied Math.}, Parthasarathy Memorial Issue, 2024, to appear.%
\newline
\noindent \lbrack BojK] R. Bojani\'{c} and J. Karamata, \textsl{On a class
of functions of regular asymptotic behavior, }Math. Research Center Tech.
Report 436, Madison, WI. 1963; reprinted in \textsl{Selected papers of Jovan
Karamata} (ed. V. Mari\'{c}, Zevod za Ud\v{z}benike, Beograd, 2009), 545-569.%
\newline
\noindent \lbrack BriD] N. Brillou\"{e}t and J. Dhombres, \'{E}quations
fonctionnelles et recherche de sous-groupes. \textsl{Aequat. Math.} \textbf{%
31} (1986), no. 2-3, 253--293.\newline
\noindent \lbrack Brz1] J. Brzd\k{e}k, Subgroups of the group Z$_{n}$ and a
generalization of the Go\l \k{a}b-Schinzel functional equation. \textsl{%
Aequat. Math. }\textbf{43} (1992), 59--71.\newline
\noindent \lbrack Brz2] J. Brzd\k{e}k, Bounded solutions of the Go\l \c{a}%
b-Schinzel equation. \textsl{Aequationes Math.} \textbf{59} (2000), no. 3,
248--254.\newline
\noindent \lbrack Chu1] J. Chudziak, Semigroup-valued solutions of the Go\l 
\k{a}b-Schinzel type functional equation. \textsl{Abh. Math. Sem. Univ.
Hamburg} \textbf{76} (2006), 91-98.\newline
\noindent \lbrack Chu2] J. Chudziak, Stability problem for the Go\l \k{a}%
b--Schinzel type functional equations. \textsl{J. Math. Anal. App.} \textbf{%
339} (2008), 454-460.\newline
\noindent \lbrack Chu3] J. Chudziak, Semigroup-valued solutions of some
composite equations. \textsl{Aequationes Math.} \textbf{88} (2014), 183--198.%
\newline
\noindent \lbrack Chu4] J. Chudziak, Continuous on rays solutions of a Go\l 
\c{a}b-Schinzel type equation. \textsl{Bull. Austral. Math. Soc.} \textbf{91}
(2015), 273--277.\newline
\noindent \lbrack Con] J. B. Conway, \textsl{A course in functional analysis}%
, Graduate texts in Math. \textbf{96}, Springer, 2$^{\text{nd.}}$ ed. 1996 (1%
$^{\text{st }}$ed. 1985).\newline
\noindent \lbrack HaaF] L. de Haan and A. Ferreira, \textsl{Extreme value
theory. An introduction. }Springer, 2006.\newline
\noindent \lbrack Jab1] E. Jab\l o\'{n}ska, Continuous on rays solutions of
an equation of the Go\l \c{a}b-Schinzel type. \textsl{J. Math. Anal. Appl.} 
\textbf{375} (2011), 223--229.\newline
\noindent \lbrack Jab2] E. Jab\l o\'{n}ska, Christensen measurability and
some functional equation. \textsl{Aequat. Math.} \textbf{81} (2011),
155--165.\newline
\noindent \lbrack Jab3] E. Jab\l o\'{n}ska, The pexiderized Go\l \c{a}%
b-Schinzel functional equation. \textsl{J. Math. Anal. Appl.} \textbf{381}
(2011), 565---572.\newline
\noindent \lbrack Jac] N. Jacobson, \textsl{Lectures in Abstract Algebra}.
vol. I. Van Nostrand, New York, 1951.\newline
\noindent \lbrack Lev] T. Levi-Civit\`{a}, Sulle funzioni che ammettono una
formula d'addizione del tipo $f(x+y)=\sum_{1}^{n}X_{j}(x)Y_{j}(y)$. \textsl{%
Atta Accad. Naz. Lincei Rend.} \textbf{22 }(1913), 181-183.\newline
\noindent \lbrack LinT] J. Lindenstrauss and L. Tzafriri, On the
complemented subspaces problem. \textsl{Israel J. Math}. \textbf{9} (1971),
263-269.\newline
\noindent \lbrack Ost1] A. J. Ostaszewski, Beurling regular variation, Bloom
dichotomy, and the Go\l \k{a}b-Schinzel functional equation. \textsl{Aequat.
Math.} \textbf{89} (2015), 725-744. \newline
\noindent \lbrack Ost2] A. J. Ostaszewski, Stable laws and Beurling kernels. 
\textsl{Adv. Appl. Probab.} \textbf{48A} (2016) (N. H. Bingham Festschrift),
239--248. \newline
\noindent \lbrack Ost3] A. J. Ostaszewski, Homomorphisms from Functional
Equations: The Goldie Equation. \textsl{Aequat. Math. }\textbf{90} (2016),
427-448 (arXiv: 1407.4089).\newline
\noindent \lbrack Ost4] A. J. Ostaszewski, Homomorphisms from Functional
Equations in Probability, in: \textsl{Developments in Functional Equations
and Related Topics}, ed. J. Brzd\k{e}k et al., Springer (2017), 171-213.%
\newline
\noindent \lbrack Pop] C. G. Popa, Sur l'\'{e}quation fonctionelle $%
f[x+yf(x)]=f(x)f(y).$ \textsl{Ann. Polon. Math.} \textbf{17} (1965), 193-198.%
\newline
\noindent \lbrack Ste] B. von Stengel, Closure properties of independence
concepts for continuous utilities. \textsl{Math. O.R.} \textbf{18} (1993),
346--389.\newline
\noindent \lbrack Stet] H. Stetkaer, \textsl{Functional Equations on groups}%
. World Scientific, 2013.\newline
\noindent \lbrack Sze] L. Sz\'{e}kelyhidi, Functional equations on abelian
groups. \textsl{Acta Math. Acad. Sci. Hungar.} 37 (1981), no. 1-3, 235--243.%
\newline

\bigskip

\noindent Mathematics Department, Imperial College, London SW7 2AZ;
n.bingham@ic.ac.uk \newline
Mathematics Department, London School of Economics, Houghton Street, London
WC2A 2AE; A.J.Ostaszewski@lse.ac.uk\newpage

\textbf{9. Appendix.}

Below, in checking routine details omitted in the main text, we offer a
fuller treatment of several results.

\bigskip

\noindent 9.1 \textit{Proof of Lemma 2.1.}

By associativity of $\circ _{\rho },$ 
\begin{eqnarray*}
K(u\circ _{\rho }v\circ _{\rho }w) &=&K(u\circ _{\rho }v)+g(u\circ _{\rho
}v)K(w) \\
&=&K(u)+g(u)K(v)+g(u\circ _{\rho }v)K(w),
\end{eqnarray*}%
and this yields $(M)$ on comparison with%
\begin{eqnarray*}
K(u\circ _{\rho }v\circ _{\rho }w) &=&K(u)+g(u)K(v\circ _{\rho }w) \\
&=&K(u)+g(u)[K(v)+g(v)K(w)] \\
&=&K(u)+g(u)K(v)+g(u)g(v)K(w).
\end{eqnarray*}%
Now%
\begin{equation*}
su\circ _{\rho }tu=su+tu+st\rho (u)u=[s+t+st\rho (u)]u,
\end{equation*}%
Put $g_{u}(t):=g(tu);$ then%
\begin{equation*}
g_{u}(s\circ _{\rho (u)}t)=g((s\circ _{\rho (u)}t)u)=g(su\circ _{\rho
}tu)=g(tu)g(tv)=g_{u}(s)g_{u}(t).
\end{equation*}%
So $g_{u}:\mathbb{G}_{\rho (u)}\mathbb{(R)\rightarrow R}_{+}$ and satisfies $%
(GFE)$ on the real line. This case is covered by established results,
e.g.[BinO5] or Th. BO below, yielding for some $\kappa (u)$%
\begin{equation*}
g_{u}(t)=(1+\rho (u)t)^{\kappa (u)/\rho (u)},
\end{equation*}%
and therefrom the cited formulas. \hfill $\square $

\bigskip

\noindent 9.2 \textit{Converse proof for Theorem 2.1.}

For $\alpha $ additive and $\beta $ constant,%
\begin{eqnarray*}
g_{\alpha ,\beta }(x\circ y) &=&e^{\alpha (x\circ y)}(1+\rho (x\circ
y))^{\beta } \\
&=&e^{\alpha (x)}e^{\alpha (y)}[(1+\rho (x))(1+\rho (y))]^{\beta } \\
&=&g_{\alpha ,\beta }(x)g_{\alpha ,\beta }(y),
\end{eqnarray*}%
so $g_{\alpha ,\beta }$ is multiplicative. \hfill $\square $

\bigskip

\bigskip

\noindent 9.3 \textit{Proof of Lemma 3.1}

We compute%
\begin{equation*}
1+\sigma \lambda (s)=g(s):\qquad \sigma =(g(1)-1)/\lambda (1)=g(1)-1:
\end{equation*}%
(i) for $\rho =0$ with $\sigma =e-1,$%
\begin{equation*}
e^{\gamma (s+t)}-1=e^{\gamma s}-1+[1+(e^{\gamma s}-1)](e^{\gamma t}-1);\text{%
\qquad\ }
\end{equation*}%
(ii) for $\rho >0,\gamma =0$ with $\sigma =\log (1+\rho )-1:$%
\begin{equation*}
\lambda (s\circ _{\rho }t)\ln (1+\rho )=\log [1+\rho (s+t+\rho st)]=\log
(1+\rho s)+\log (1+\rho t)\text{;}
\end{equation*}%
(iii) for $\rho >0$ with $\sigma :=(1+\rho )^{\gamma /\rho }-1$%
\begin{equation*}
\lambda (s\circ _{\rho }t)=\lambda (s)\circ _{\sigma }\lambda (t)\text{ }.
\end{equation*}

\hfill $\square $

\bigskip

\noindent 9.4 \textit{Direct proof of Theorem 3.1 assuming Gateaux
differentiability.}

We begin as in Lemma 4.1. As there, for fixed $u\in X,$ put%
\begin{equation*}
x\circ _{h(u)}y:=x+h(u)y,\qquad a\circ _{g(u)}b:=a+g(u)b.
\end{equation*}%
Context permitting, abbreviate this to $\circ _{g}$etc. Starting from $u$
and $v:=K(u),$ we again define a pair of sequences of `powers', by iterating
the two operations $\circ _{h}$ and $\circ _{g}$. These iterates are defined
inductively: 
\begin{equation*}
u_{h}^{n+1}=u\circ _{h}u_{\rho }^{n},\qquad v_{g}^{n+1}=v\circ _{g}v_{\rho
}^{n},\qquad \text{with }u_{\rho }^{1}=u,\qquad v_{g}^{1}=v.
\end{equation*}%
Then, for $n\geq 1,$%
\begin{equation}
K(u_{h}^{n+1})=K(u)+g(u)K(u_{h}^{n})=K(u)_{g}^{n+1}=v_{g}^{n+1}. 
\tag{$\ast
$}
\end{equation}%
Motivated by the case%
\begin{equation*}
K(u_{h}^{2})=K([1+h(u)]u)=K(u)+g(u)K(u)=[1+g(u)]K(u),
\end{equation*}%
the recurrence $(\ast )$ justifies associating with the iterates above a
sequences of `coefficients' $(g_{n}(.)),$ $(h_{n}(.))$, by writing 
\begin{equation*}
v_{g}^{n}=g_{n}(u)K(u),\qquad u_{h}^{n}=h_{n}(u)u:\qquad
K(h_{n}(u))=g_{n}(u)K(u).
\end{equation*}%
Solving appropriate recurrences arising from $(\ast )$ above for the
iterations $u_{h}^{n+1}=u\circ _{h}u_{\rho }^{n}$ and $v_{g}^{n+1}=v\circ
_{g}v_{\rho }^{n}$ gives%
\begin{eqnarray*}
u_{h}^{n} &=&h_{n}(u)u=\frac{h(u)^{n}-1}{h(u)-1}u\text{ or }nu\text{ if }%
h(u)=1, \\
v_{g}^{n} &=&g_{n}(u)K(u)=\frac{g(u)^{n}-1}{g(u)-1}K(u)\text{ or }nK(u)\text{
if }g(u)=1.
\end{eqnarray*}%
Replacing $u$ by $u/n,$%
\begin{equation*}
K(u/n)=g_{n}(u/n)^{-1}K(h_{n}(u/n)u/n),
\end{equation*}%
so%
\begin{equation}
K(h_{m}(u/n)u/n)=g_{m}(u/n)K(u/n)=g_{m}(u/n)g_{n}(u/n)^{-1}K(h_{n}(u/n)u/n).
\tag{$\ast \ast $}
\end{equation}%
Here%
\begin{equation*}
g_{m}(u/n)g_{n}(u/n)^{-1}=(g(u/n)^{m}-1)/(g(u/n)^{n}-1)\text{ or }m/n\text{
if }g(u/n)=1.
\end{equation*}

First suppose for all large $n$ that $g(u/n)\neq 1.$ Fix $t.$ For each $n\in 
\mathbb{N}$ choose $m=m(n)\in \mathbb{N},$ so that $t(n):=m(n)/n\rightarrow
t,$ and write%
\begin{equation*}
g(u/n)-g(0)=g_{u}^{\prime }(\varepsilon _{n})(1/n),
\end{equation*}%
by the mean-value theorem with $0<\varepsilon _{n}<1/n.$ Then as $%
n\rightarrow \infty ,$ since $g_{u}^{\prime }(\varepsilon _{n})\rightarrow
g_{u}^{\prime }(0),$%
\begin{eqnarray*}
g(u/n)^{t(n)n} &=&[1+g^{\prime }(\varepsilon _{n})(u/n)]^{t(n)n} \\
&=&\exp \{(t+\delta _{n})u[\log [1+g_{u}^{\prime }(\varepsilon
_{n})(u/n)]/(u/n)] \\
&\rightarrow &\exp tg^{\prime }(0).
\end{eqnarray*}%
So 
\begin{eqnarray*}
g_{m}(u/n)/g_{n}(u/n) &=&\frac{g(u/n)^{t(n)n}-1}{g(u/n)^{n}-1}, \\
&=&\frac{\exp \{t(n)[n\log [1+g_{u}^{\prime }(\varepsilon _{n})(1/n)]\}-1}{%
\exp \{t(n)[n\log [1+g_{u}^{\prime }(u\varepsilon _{n})(1/n)]\}-1} \\
&\rightarrow &[e^{g_{u}^{\prime }(0)t}-1]/[e^{g_{u}^{\prime }(0)}-1].
\end{eqnarray*}%
Note that this still holds even if $g(u/n)=1$ for infinitely many $n.$ For
then $g_{m}(u/n)=n,$ and $g_{m}(u/n)/g_{m}(u/n)=m/n\rightarrow t;$
furthermore,%
\begin{equation*}
0=n[g(u/n)-1]=g_{u}^{\prime }(\varepsilon _{n})\rightarrow g_{u}^{\prime
}(0),
\end{equation*}%
i.e. $g_{u}^{\prime }(0)=0$. Likewise, writing 
\begin{equation*}
\lbrack h(u/n)-1]=h_{u}^{\prime }(\varepsilon _{n})(1/n),\qquad
n[h(u/n)-1]=h_{u}^{\prime }(\varepsilon _{n})
\end{equation*}%
\begin{eqnarray*}
h_{m}(u/n)/n &=&\frac{h(u/n)^{m}-1}{h(u/n)-1}\frac{1}{n}, \\
&=&\frac{\exp \{t(n)[n\log [1+h_{u}^{\prime }(u\varepsilon _{n})(1/n)]\}-1}{%
h^{\prime }(0)}, \\
&\rightarrow &[e^{h_{u}^{\prime }(0)t}-1]/h_{u}^{\prime }(0).
\end{eqnarray*}

Again this still holds even if $h(u/n)=1$ for infinitely many $n.$ For then $%
h_{m}(u/n)=m,$ and $h_{m}(u/n)/n=m/n\rightarrow t;$ furthermore, again $%
h_{u}^{\prime }(0)=0.$%
\begin{equation*}
K\left( \frac{e^{\rho (u)t}-1}{\rho (u)}u\right) =\frac{e^{\gamma (u)t}-1}{%
e^{\gamma (u)}-1}K(u),
\end{equation*}%
Solving for $s$ and substituting for $t$ in terms of $s$ gives for $%
c(u):=\gamma (u)/\rho (u)$ 
\begin{eqnarray*}
s &=&\frac{e^{\rho (u)t}-1}{\rho (u)}\qquad t(s):=\log [1+s\rho (u)]/\rho
(u)\qquad \\
\lambda _{u}(s) &:&=\frac{e^{\gamma (u)t}-1}{e^{\gamma (u)}-1}=\frac{%
e^{c(u)\log [1+s\rho (u)]}-1}{e^{c(u)\rho (u)}-1}=\frac{(1+s\rho
(u))^{c(u)}-1}{(1+\rho (u))^{c(u)}-1}.
\end{eqnarray*}%
This yields%
\begin{equation*}
K(su)=\lambda _{u}(s)K(u).
\end{equation*}

\begin{eqnarray*}
K(su+h(su)tu) &=&K(su)+g(su)K(tu): \\
\lambda _{u}(s+h(su)t)K(u) &=&[\lambda _{u}(s)+g(su)\lambda _{u}(t)]K(u)
\end{eqnarray*}%
The final assertion is now immediate. \hfill $\square $

\bigskip

\noindent 9.5 \textit{Proof of Lemma 7.1.}

For given $K,$ suppose both $(K,g)$ and $(K,h)$ solve $(GFE)$. Fix $v\in X$
with $K(v)\neq 0$; then for arbitrary $u\in X$%
\begin{equation*}
K(u)+h(u)K(w)=K(u\circ _{\rho }v)=K(u)+g(u)K(w),
\end{equation*}%
so $g(u)=h(u).$ So if $g^{\tau }=g^{\sigma },$ then 
\begin{equation*}
1+\sigma (K(u))=g^{\sigma }(u)=g^{\tau }(u)=1+\tau (K(u)).
\end{equation*}%
Furthermore,%
\begin{eqnarray*}
\sigma (K(u\circ _{\rho }v)) &=&g(u)g(v)-1=(1+\sigma (K(u))(1+\sigma (K(v))-1
\\
&=&\sigma (K(u))+[1+\sigma (K(v))]\sigma (K(v)) \\
&=&\sigma (K(u))\circ _{1}\sigma (K(v))\qquad \qquad \qquad \square
\end{eqnarray*}

\bigskip

\noindent 9.6 \textit{Alternative approach to the Tetrachotomy of Theorem
8.1.}

One can arrive at the cases given in Theorem 8.1 by a more transparent,
albeit pedestrian, route, assuming either Gateaux differentiability or,
equivalently, finiteness of the limits $\gamma (u)=\lim_{n}n\delta _{n}^{g}$
and $\rho (u)=\lim_{n}n\delta _{n}^{h}.$

\bigskip

\noindent \textbf{Proposition. }\textit{For }$K,h,g$ \textit{continuous and} 
$(K,h,g)$ \textit{satisfying }$(GGE)$ \textit{with }$K$ \textit{non-trivial
(in the sense that there is }$w\in X$ \textit{with }$K(w)\neq 0$ \textit{and
also }$\lambda _{w}\neq \mathrm{id)}$ \textit{and for any }$u\in X$ \textit{%
with }$K(u)\neq 0,$\textit{\ one of the following cases occurs.}

\noindent (i) \textit{If} $\rho (u)=\gamma (u)=0,$ \textit{then} 
\begin{equation*}
g_{u}(s)=h_{u}(s)=1;
\end{equation*}

\noindent (ii) \textit{If} $\rho (u)=0$\textit{\ and} $\gamma (u)\neq 0,$%
\textit{\ then }%
\begin{equation*}
h(su)=1\text{ and }g(su)=e^{s\gamma (u)};
\end{equation*}

\noindent (iii) \textit{If} $\rho (u)\neq 0$\textit{\ and }$\gamma (u)=0,$%
\textit{\ then}%
\begin{equation*}
h_{u}(s)=1+s\rho (u),\qquad g_{u}(s)=1;
\end{equation*}

\noindent (iv) \textit{If} $\gamma (u)\neq 0\neq \rho (u),$\textit{\ then}%
\begin{equation*}
h_{u}(u)=(1+s\rho (u)),\qquad g_{u}(s)=(1+s\rho (u))^{(\gamma /\rho )}.
\end{equation*}

\bigskip

\noindent \textbf{Proof. }Fix $u$ with $K(u)\neq 0.$ We consider the four
cases identified in Theorem 8.1 which $u$ gives rise to. Consider $s,t\in 
\mathbb{R}$ and set\textbf{\ }$v=tu.$ By Theorem 3.1 from $K(su)\neq 0,$
enabling its cancellation,%
\begin{eqnarray*}
K(su+th(su)su) &=&K([1+th(su)]su): \\
K(su+th(su)su) &=&K(su)+g(su)K(tsu), \\
\lambda _{su}(1+th(su))K(su) &=&K(su)+g(su)\lambda _{su}(t)K(su): \\
\lambda _{su}(1+th(su)) &=&1+g(su)\lambda _{su}(t).
\end{eqnarray*}%
Differentiating w.r.t. $t$ one obtains for general $t$ and also for $t=0:$%
\begin{eqnarray*}
\lambda _{su}^{\prime }(1+th(su))h(su) &=&g(su)\lambda _{su}^{\prime }(t), \\
\lambda _{su}^{\prime }(1)h(su) &=&g(su)\lambda _{su}^{\prime }(0).
\end{eqnarray*}%
We reduce $\lambda _{su}^{\prime }$ to $\lambda _{u}^{\prime }$ by recalling
from Corollary 4.3 that%
\begin{equation*}
\lambda _{su}^{\prime }(t)=\lambda _{u}^{\prime }(st)s/\lambda _{u}(s).
\end{equation*}%
This transforms the preceeding two equations to%
\begin{equation}
\lambda _{u}^{\prime }(s+sth(su))h(su)=g(su)\lambda _{u}^{\prime }(st), 
\tag{$D1$}
\end{equation}%
and%
\begin{equation}
\lambda _{u}^{\prime }(s)h(su)=g(su)\lambda _{u}^{\prime }(0)  \tag{$D2$}
\end{equation}%
We deduce the form of $g$ and $h$ from the above two equations $(D1,2)$ by
computing $\lambda _{u}^{\prime }(s)$ case by case according to the four
defining clauses of $\lambda _{u}(t)=\lambda (t;\rho (u),\gamma (u)).$

\bigskip

\noindent \textit{Case }(i). $\lambda (t)\equiv t$ and $\rho (u)=\gamma
(u)=0,$ for which $\lambda _{u}^{\prime }=1.$ Equation $(D2)$ reduces to 
\begin{equation*}
h(su)=g(su)\text{.}
\end{equation*}%
By Proposition 4.4, $g$ satisfies $(GS)$ and so for some $\delta =\delta (u)$
\begin{equation*}
h_{u}(s)=g_{u}(s)=1+s\delta _{u}.
\end{equation*}%
Here $\delta _{u}=0,$ since%
\begin{equation*}
\delta _{u}=h_{u}^{\prime }(0)=g_{u}^{\prime }(0)=\rho (u)=\gamma (u)=0.
\end{equation*}

\noindent \textit{Case} (ii). $\lambda _{u}(t)=(e^{\gamma t}-1)/(e-1)$ and $%
\rho (u)=0,\gamma (u)\neq 0\gamma (u)\neq 0,$ for which%
\begin{equation*}
\lambda _{u}^{\prime }(s)=\gamma e^{\gamma s}/(e-1)\text{ and }\lambda
_{u}^{\prime }(0)=\gamma /(e-1).
\end{equation*}%
Here, since $\gamma (u)\neq 0$, equation $(D2)$ reduces to%
\begin{equation*}
\gamma e^{\gamma s}h(su)=g(su)\gamma :\qquad e^{\gamma s}h(su)=g(su).
\end{equation*}%
Then substituting for $g$ into $(D1)$ gives $h_{u}\equiv 1,$ since%
\begin{equation*}
\lambda _{u}^{\prime }(s)=\gamma e^{\gamma s}/(e-1)
\end{equation*}

Here equation $(D2)$ reduces to%
\begin{eqnarray*}
\gamma (su)e^{\gamma (su)}h(su) &=&\gamma (su)g(su),\qquad \text{(cancel by }%
\gamma (u)\neq 0\text{)} \\
e^{s\gamma (u)}h(su) &=&g(su).
\end{eqnarray*}

Then substituting for $g$ into $(D1)$ gives $h_{u}\equiv 1.$ Indeed, $(D1)$
gives%
\begin{eqnarray*}
\frac{\lambda _{u}^{\prime }(s+sth(su))s}{\lambda _{u}(s)}h(su) &=&g(su)%
\frac{\lambda _{u}^{\prime }(st)s}{\lambda _{u}(s)}: \\
\lambda _{u}^{\prime }(s+sth(su))h(su) &=&g(su)\lambda _{u}^{\prime }(st),
\end{eqnarray*}%
and substitution for $g(su)$ and cancelling $h(su)$ yields%
\begin{equation*}
\lambda _{u}^{\prime }(s+sth(su))=e^{s\gamma (u)}\lambda _{u}^{\prime }(st).
\end{equation*}%
Substitution for $\lambda _{u}^{\prime }$ and cancelling the non-zero $%
\gamma $ and $(e-1)^{-1}$ gives 
\begin{eqnarray*}
e^{\gamma \lbrack s+sth(su)]} &=&e^{s\gamma (u)}e^{\gamma st} \\
\gamma sth(su) &=&\gamma st,\qquad \text{ i.e. }h(su)=1.
\end{eqnarray*}

\noindent \textit{Case} (iii). $\lambda _{u}(t)=\ln (1+\rho (u)t)/\ln
(1+\rho (u))$ and $\gamma (u)=0,$ for which $\lambda _{u}^{\prime }(t)=\rho
(u)/[(1+\rho (u)t)\ln (1+\rho (u))].$ Here $D(2)$ gives%
\begin{equation*}
\frac{h(su)}{1+s\rho (u)}=g(su)\text{.}
\end{equation*}%
From here substitution for $g$ into $(D1)$ and cancelling $h$ gives%
\begin{eqnarray*}
\frac{1}{1+\rho (u)[s+sth(su)]} &=&\frac{1}{1+st\rho (u)}\frac{1}{1+s\rho (u)%
},\qquad \text{(reciprocate:)} \\
1+s\rho (u)[1+th(su)] &=&[1+s\rho (u)][1+st\rho (u)], \\
1+s\rho (u)+st\rho (u)h(su) &=&1+s\rho (u)+st\rho (u)+s^{2}t\rho (u)^{2}, \\
st\rho (u)h(su) &=&st\rho (u)[1+s\rho (u)]:\qquad \lbrack 1+\rho (u)s]=h(su).
\end{eqnarray*}

\noindent \textit{Case} (iv). $\lambda _{u}(t)=[(1+t\rho )^{\gamma /\rho
}-1]/[(1+\rho )^{\gamma /\rho }-1]$ with $\gamma (u)\neq 0\neq \rho (u),$
for which%
\begin{equation*}
\lambda _{u}^{\prime }(t)=[\gamma (1+t\rho )^{\gamma /\rho -1}]/[(1+\rho
)^{\gamma /\rho }-1]
\end{equation*}%
Here, $(D1)$ gives that%
\begin{equation*}
\lambda _{u}^{\prime }(s+sth(su))h(su)=g(su)\lambda _{u}^{\prime }(st),
\end{equation*}%
\begin{equation*}
h(su)\gamma (u)(1+\rho (u)s[1+th(su)])^{(\gamma /\rho )-1}=g(su)\gamma
(u)(1+st\rho (u))^{(\gamma /\rho )-1},
\end{equation*}%
and $(D2)$ that%
\begin{equation*}
h(su)(1+s\rho (u))^{(\gamma /\rho )-1}=g(su).
\end{equation*}%
From here substitution for $g(su)$ in ($D1)$ with $s$ fixed, gives 
\begin{equation*}
h(su)(1+\rho (u)s[1+th(su)])^{(\gamma /\rho )-1}=h(su)(1+s\rho (u))^{(\gamma
/\rho )-1}(1+st\rho (u))^{(\gamma /\rho )-1}.
\end{equation*}%
Cancelling $h(su),$ and using the identity 
\begin{equation*}
\lbrack (1+s\rho )(1+st\rho )]=1+s\rho \lbrack 1+t+ts\rho ],
\end{equation*}%
\begin{eqnarray*}
(1+\rho (u)s[1+th(su)])^{(\gamma /\rho )-1} &=&[(1+s\rho (u))(1+ts\rho
(u))]^{(\gamma /\rho )-1}, \\
&=&[1+s\rho (u)[1+t+ts\rho (u))]]^{(\gamma /\rho )-1}.
\end{eqnarray*}%
Disposal of the exponent yields%
\begin{eqnarray*}
1+\rho (u)s[1+th(su)] &=&1+s\rho (u)[1+t(1+s\rho (u))], \\
h(su) &=&1+s\rho (u).
\end{eqnarray*}%
But%
\begin{equation*}
h(su)(1+s\rho (u))^{(\gamma /\rho )-1}=g(su),
\end{equation*}%
completing the claim. \hfill $\square $

\end{document}